\documentclass[11pt]{article}

\usepackage{latexsym}
\usepackage{amsmath, amsfonts,amssymb, amsthm, euscript,makeidx,color,mathrsfs}
\usepackage{comment}

\newtheorem{assumption}{Assumption}
\def\qed{ \ \vrule width.2cm height.2cm depth0cm\smallskip}

\newcommand{\hP}{\hat\dbP}

\newcommand{\ba}{\begin{array}}
\newcommand{\ea}{\end{array}}
\newcommand{\be}{\begin{equation}}
\newcommand{\ee}{\end{equation}}
\newcommand{\bea}{\begin{eqnarray}}
\newcommand{\eea}{\end{eqnarray}}
\newcommand{\beaa}{\begin{eqnarray*}}
\newcommand{\eeaa}{\end{eqnarray*}}

\def\neg{\negthinspace}

\def\dbE{\mathbb{E}}
\def\dbF{\mathbb{F}}

\def\dbL{\mathbb{L}}

\def\dbN{\mathbb{N}}
\def\dbP{\mathbb{P}}
\def\dbR{\mathbb{R}}

\def\a{\alpha}
\def\b{\beta}
\def\g{\gamma}
\def\d{\delta}
\def\e{\varepsilon}
\def\z{\zeta}
\def\k{\kappa}
\def\l{\lambda}

\def\si{\sigma}
\def\t{\tau}
\def\f{\varphi}
\def\th{\theta}
\def\o{\omega}
\def\h{\widehat}
\def\D{\Delta}
\def\Th{\Theta}
\def\L{\Lambda}

\def\O{\Omega}
\def\cD{{\cal D}}

\def\cF{{\cal F}}

\def\cH{{\cal H}}
\def\cI{{\cal I}}

\def\cO{{\cal O}}

\def\hE{\mathbb{E}}
\def\hF{\mathbb{F}}

\def\hN{\mathbb{N}}

\def\hP{\mathbb{P}}

\def\hR{\mathbb{R}}

\def\no{\noindent}

\def\ms{\medskip}
\def\bs{\bigskip}
\def\q{\quad}
\def\qq{\qquad}

\def\pa{\partial}
\def\cd{\cdot}
\def\cds{\cdots}
\def\lan{\langle}
\def\ran{\rangle}

\def\bz{{\bf z}}

\def\tr{\hbox{\rm tr}}

\def\qed{ \hfill \vrule width.25cm height.25cm depth0cm\smallskip}

\newcommand{\basa}{\begin{assumption}}
\newcommand{\easa}{\end{assumption}}

\newcommand{\bas}{\begin{assum}}
\newcommand{\eas}{\end{assum}}

\def\ua{\mathop{\uparrow}}
\def\da{\mathop{\downarrow}}

\def\lan{\mathop{\langle}}
\def\ran{\mathop{\rangle}}

\def\pa{\partial}
\def\h{\widehat}

 \def\cd{\cdot}
\def\cds{\cdots}

\def\tr{\hbox{\rm tr$\,$}}

\def\dis{\displaystyle}

\def\cad{c\`{a}dl\`{a}g}

\def\1{{\bf 1}}

\def\:{\!:\!}
\def\reff#1{{\rm(\ref{#1})}}

at 9pt

\begin{document}

\newtheorem{thm}{Theorem}[section]
\newtheorem{lem}[thm]{Lemma}
\newtheorem{cor}[thm]{Corollary}
\newtheorem{prop}[thm]{Proposition}
\newtheorem{rem}[thm]{Remark}
\newtheorem{eg}[thm]{Example}
\newtheorem{defn}[thm]{Definition}
\newtheorem{assum}[thm]{Assumption}

\renewcommand {\theequation}{\arabic{section}.\arabic{equation}}
\def\thesection{\arabic{section}}

\title{\bf  Pathwise Taylor Expansions for Random Fields on Multiple Dimensional Paths}

\author{
Rainer Buckdahn\thanks{ \no D\'{e}partement de Math\'{e}matiques,
Universit\'{e} de Bretagne-Occidentale, F-29285 Brest Cedex, France; and School of Mathematics, Shandong University, Jinan 250100, 
P.R. China. Email: Rainer.Buckdahn@univ-brest.fr. },  ~ Jin Ma\thanks{ \noindent Department of
Mathematics, University of Southern California, Los Angeles, 90089, USA.
%Department of Mathematics, Purdue University, West Lafayette, IN
%47907-1395;
Email: jinma@usc.edu. This author is supported in part
by US NSF grant \#1106853. } ~ and ~{Jianfeng Zhang}\thanks{\noindent
Department of Mathematics, University of Southern California, Los
Angeles, CA 90089, USA. E-mail: jianfenz@usc.edu. This author is
supported in part by NSF grant \#1008873. }}
% Part of this work was
%completed while the second and third author were visiting the School of
%Mathematics and QiLu Institute of Finance, Shandong University, China, whose
%hospitality was greatly appreciated.}}}

\date{}
%\date{\today}
\maketitle

\begin{abstract}
In this paper we establish the pathwise Taylor expansions for random fields that are ``regular"
 in the spirit of Dupire's path-derivatives \cite{Dupire}. Our result is motivated by but extends  the recent
result  of Buckdahn-Bulla-Ma \cite{BBM}, when translated into the language of pathwise calculus. We show that with such a language 
the pathwise Taylor expansion can be naturally carried out to any order and for any dimension, and it coincides with the existing results 
when  reduced to these special settings. 
% the appears exactly in the same form as the standard (deterministic) Taylor expansion, and .
% ``without tear". 
More importantly, the expansion can be both ``forward" and "backward" (i.e., the temporal increments can be both positive and negative), and the remainder is estimated in a pathwise manner. This result will be the main building block for our new notion of viscosity solution to forward path-dependent PDEs corresponding to (forward) stochastic PDEs in our accompanying paper \cite{BMZ}.

\end{abstract}

\vfill \bs

\no

{\bf Keywords.} \rm Path derivatives, pathwise Taylor expansion, functional It\^{o} formula, It\^{o}-Ventzell formula, stochastic partial differential equations.
%  stochastic viscosity solutions,
%Doss transformation,
%stochastic characteristics.

\bs

\no{\it 2000 AMS Mathematics subject classification:} 60H07,15,30;
35R60, 34F05.

\eject

\section{Introduction}
\label{sect-Introduction}
\setcounter{equation}{0}

In this paper we are interested in establishing the pathwise Taylor expansions for the It\^{o}-type random field of the form
 \bea
 \label{itofield}
 u(t,x)=u_0(x)+\int_0^t
 \a(s,x)ds+\int_0^t  \b(s,x) \circ dB_s,
  \eea
where $B$ is a $d$-dimensional standard Brownian motion, defined on a
complete probability space $(\O,\cF,\hP)$, and ``$\circ$" means Stratonovic integral. In particular we are interested in
 such expansions for the solution to the 
following fully nonlinear stochastic partial differential equations (SPDE):
 \bea
 \label{SPDEintro}
 u(t,x)=u_0(x)+\int_0^t
 f(s,x,\cd, u, \pa_xu, \pa_{xx}u)ds+\int_0^t g(s,x,\cd, u, \pa_xu) \circ dB_s,
 \eea
where $f$ and $g$ are  random fields that have certain regularity in their spatial variables.

In our previous work \cite{BM3} we studied the so-called
{\it pathwise stochastic Taylor expansion} for a class of
It\^{o}-type random fields. The main result can be briefly described
as follows. 
 Suppose that $u$ is
a random field of the form (\ref{itofield}), and $B$ is a one dimensional Brownian motion. If we denote
$\hF=\{\cF_t\}_{t\ge0}$ to be the natural filtration generated by
$B$ and augmented by all $\hP$-null sets in $\cF$,  then under
reasonable regularity assumptions on the integrands $\a$ and $\b$,
 the following stochastic ``Taylor expansion" holds:
{\it For any stopping time $\t$ and any $\cF_\t$-measurable,
square-integrable random variable $\xi$, and for any sequence of
random variables $\{(\t_k,\xi_k)\}$ where $\t_k$'s are stopping
times such that either $\t_k>\t$, $\t_k\da \t$; or $\t_k<\t$,
$\t_k\ua \t$,
%such that $\cF_{\t_k}\subseteq\cF_\t$
and $\xi_k$'s are all $\cF_{\t_k\wedge \t}$-measurable, square
integrable random variables, converging to $\xi$ in $L^2$, it holds
%almost surely 
that
 \bea
 \label{itotaylor}
 &&u(\t_k,\xi_k) =u(\t,\xi)+a(\t_k-\t)+b(B_{\t_k}-B_\t)+p(\xi_k-\xi) \\
 &&\q +{c\over
 2}(B_{\t_k} -B_\t)^2+ q(\xi_k-\xi)(B_{\t_k}-B_\t)+\frac12  X(\xi_k-\xi)^2 + o(|\t_k-\t|+|\xi_k-\xi|^2), \nonumber
 \eea
where $(a,b,c,p,q,X)$ are all $\cF_\t$-measurable random variables,
and the remainder $o(\z_k)$ are such that $o(\z_k)/{\z_k} \to 0$ as
$k\to\infty$, in probability.} Furthermore, the six-tuple
$(a,b,c,p,q,X)$ can be determined explicitly in terms of $\a$,
$\b$ and their derivatives (in certain sense).

While the Taylor expansion (\ref{itotaylor}) reveals the possibility of estimating the remainder in
a stronger form than mean-square (cf. e.g., \cite{KloPla}), it is not satisfactory for the study of pathwise property of the random fields
which is essential in the study of, e.g., stochastic viscosity solution. In a subsequent paper 
(Buckdahn-Bulla-Ma \cite{BBM}) the result was extended to the case where the expansion
could be made around any random time-space point $(\t, \xi)$ where $\t$ does not have to be a
stopping time; and more importantly, the remainder was estimated in a pathwise manner, in the spirit of the
Kolmogorov continuity criterion. In other words, modulo a $\hP$-null set, the estimate holds for each
$\o$, locally uniformly in $(t,x)$. Furthermore, all the coefficients can be calculated explicitly in terms
of a certain kind of ``derivatives" for It\^o-type random field, introduced in \cite{BBM} (see more detailed description
in \S8 of this paper). It is noted, however,
that a main drawback of the result in \cite{BBM} is that the derivatives involved are not intuitive, and are difficult to verify.
A more significant weakness of the result is that the dimension of the Brownian motion is restricted to 1, which, as we 
shall see in this paper, reduced the complexity of the Taylor expansion drastically.
%which greatly limited its applicability. 

The main purpose of this paper is to re-investigate the Taylor expansion in a much more general setting, but with 
a different ``language." 
In particular, we shall allow both the spatial variable and the Brownian motion to be multi-dimensional, and the random field is ``regular" in a very different way. To be more precise, we shall introduce
a new notion of ``path-derivative" in the spirit of Dupire \cite{Dupire}  to  
% between this paper and our previous ones is that we shall
impose a different type of regularity, that is, the regularity on the variable $\o\in\O$. 
%he spirit of Dupire's derivatives \cite{Dupire}.  
Such a language turns out to be very effective, and many originally cumbersome expressions in stochastic
analysis becomes intuitive and very easy to understand. For example, even without using the Stratonovic
integral, the It\^o-Ventzell formula reads exactly like the multi-dimensional It\^o formula, and both integrands
for the Lebesgue integral and the stochastic integral can be memorized simply as ``chain rule", with respect to 
time and path, respectively. For this reason we shall name it ``pathwise It\^o-Ventzell formula" (see Section 4 
below for details). We should note, however, that our  path derivative is much weaker than the original
one by Dupire (see also \cite{CF}), and applies  to all semi-martingales. But on the other hand such a 
generality brings out some intrinsic ``rough-path" nature of the Brownian motion. Among other things, for example, 
the ``Hessian" under the current path-derivatives will be {\it asymmetric} in a  general multi-dimensional setting, reflecting 
the nature of {\it L\'evy area} in the rough-path theory (cf. e.g., \cite{Lyons}, or  \cite{FV}). 

We would like to point out that the Taylor expansion for stochastic processes, especially for the solutions of stochastic
differential equations, is not new. There is a  large amount of literature on the subject, from various perspectives. 
We refer to the books of Kloeden-Platten \cite{KloPla} from the numerical approximation point of view; and of 
Friz-Victoir \cite{FV} from the rough-path point of view, as well as the numerous references cited therein. In fact,
all Taylor expansions resemble each other in their forms, language notwithstanding, and the difference often lies in the 
error (remainder) estimates. The main feature of our results is the following.
%in this paper and the existing ones. 
First, our Taylor expansion applies to general random fields
and stochastic processes, and therefore does not depend on the special structure for being a solution to a differential
equation, whence ``non-Markovian" in nature. Second, unlike our previous work, we shall provide a unified treatment
of the Taylor expansion up to any order, and allowing the temporal increment to be both 
``forward" and ``backward". Finally, and most importantly, we pursue the pathwise estimate for the
remainder, that is, the error of the expansion is estimated uniformly for all paths $\o$, modulo a common null set. 
The main difficulty, compared to an $L^2$ estimate (or  in the 
sense of ``in probability") that we often see in the literature, is that one cannot use the isometry between the 
$L^2$-norms of stochastic integrals and the  $L^2$-norms of the Lesbesgue integrals, thus it requires some
novel treatments of multiple integrals.
%, and the results will be truly pathwise. 
The trade-off for being able to do this, however,  is that we require some new regularities of the random field with
respect to the ``paths". These requirements, when unified under our new language of ``pathwise analysis",  are direct and 
easy to check. 
%the problems becomes more straight forward. 
To our best knowledge, the pathwise Taylor expansion in such a generality is new. 

%of the nature of the path behavior.rewe will see
%some very interesting the  and all the technical assumptions
%become very straightforward and natural, and more importantly, the pathwise Taylor
%expansion appears exactly in the same form as the standard (deterministic) Taylor expansion(!). 
%Moreover, we shall check that in the case when all the coefficients of the SPDE (\ref{SPDEintro}),
%the more or less ``mysterious" coefficients in the recent result 
%of  \cite{BBM}, which motivated the present work, exactly coincides with our new result when translated into the language of 
%pathwise calculus. Thus the result of this paper is indeed a non-trivial extension of our previous ones.  

It is worth noting that our Taylor expansion is the first step of our  study of the  viscosity solution to the (forward) ``path-dependent PDEs"
(PPDEs) corresponding to the forward SPDE (\ref{SPDEintro}), which will be the main topic of our accompanying paper \cite{BMZ}.
We would only like to comment here that a classical solution in the traditional sense does not necessarily permit a pathwise Taylor
expansion. Therefore a somewhat  convoluted treatment of the solution to the SPDEs will have to be carried out based on the pathwise
Taylor expansions, as it was seen in the deterministic viscosity solution theory as well as the existing studies of stochastic viscosity solutions 
(cf. e.g., \cite{BM3}).

The rest of the paper is organized as follows. In Section 2 we give
all necessary notations and
%. In Section \ref{Notion of stochastic viscosity solution} we motivate
%and give the definition of the stochastic viscosity solution. For these reasons, we
introduce the definition of  ``path-derivatives." In Section
3 we give a heuristic analysis for a simpler case, the second order  expansion for It\^o processes,
to illustrate the main points of our method. In  section 4 we prove the crucial estimates for the remainders of higher order Taylor expansions.
%estimate that is crucial to the main result. 
In section 5 we extend  the Taylor expansion
to It\^o random fields; and in Section 6 we weaken the regularity assumptions of the coefficients to 
H\"older spaces. % we prove the main results,
In Section 7 we apply the Taylor expansion to the solutions to stochastic
PDEs, and  finally, in Section 8 we compare the main theorem with our previous result \cite{BBM}. 
%In section 5 we present a general version for higher order expansions. 
%Finally, in sections 7 we discuss some regularity issues for the classical 
%solutions of SPDE (\ref{SPDEintro}), which more or less justifies the assumptions
%of the main theorem.
%

\section{Preliminaries}
\label{sect-Temp}
\setcounter{equation}{0}

Throughout this paper we denote $\O:= \big\{\o\in C([0,\infty), \dbR^d): \o_0={\bf 0}\big\}$ to be
the set of continuous paths starting from the origin,  $B$ the canonical process on $\O$, 
$\dbP_0$ the Wiener measure,  $\dbF = \{\cF_t\}_{t\ge 0}$ the $\dbP_0$-augmented filtration generated by $B$, and $\L := [0,\infty)\times \O$. Here and in the sequel, we use ${\bf 0}$ to denote vectors or matrices with appropriate dimensions whose components are all equal to $0$, and for any dimension $m$, we take the convention that $\hR^m := \hR^{m\times 1}$ denotes the set of column vectors.  Define
 \beaa
 x \cd x' := \sum_{i=1}^m x_i x'_i 
 ~~\mbox{for any}~~x, x' \in \dbR^m,
 ~~\g : \g' := \tr[\g (\g')^T]
 ~~\mbox{for any}~~\g, \g'\in \dbR^{m\times n},
 \eeaa
and $|x|^2 := x\cd x$, $|\g|^2 := \g : \g$. Here $^T$ denotes the transpose.

\subsection{ Path derivatives for It\^{o} processes}
% and functional It\^o formula}

Let $\dbL^0(\L, \dbR^{m\times n})$ denote the set of $\dbF$-progressively measurable  processes $u: \L \to \dbR^{m\times n}$, and $\dbL^0(\L) := \dbL^0(\L, \dbR)$.
 Strongly motivated by the functional It\^{o} formula initiated by Dupire \cite{Dupire} (see also Cont and Fournie
 \cite{CF} and a slight variation by Ekren-Touzi-Zhang \cite{ETZ}), in what follows we introduce the notion of ``path-derivatives", which will
be the foundation of our
pathwise stochastic analysis.

Recall that $u\in\dbL^0(\L)$ is a  semimartingale if there exist $A\in \dbL^0(\L)$ and $\beta\in \dbL^0(\L, \dbR^d)$ such that  
\bea
\label{uito} 
 u_t=u_0+A_t+\int_0^t \b_s \cd dB_s, ~\mbox{and}~ V_0^t(A) + \int_0^t |\b_s|^2 ds <\infty,  \q  t\ge 0, ~\hP_0\mbox{-a.s.,}
 \eea
where $V_0^t(A)$ is the total variation of $A$ on $[0, t]$.

%We now define the {\it generalized path-derivative} of $u$  in $(t,\o)$ as follows.
\begin{defn}
\label{Dw}
Let $u$ be a  semimartingale in the form of  (\ref{uito}). We define:
\bea
\label{Dou}
\pa_\o u:= \beta.
\eea
Moreover, if $\beta$ is  a  semimartingale and $A_t = \int_0^t \a_s ds$ for some $\a\in \dbL^0(\L)$ , then we define
\bea
\label{Dtu}
\pa^2_{\o\o}u := \pa_\o \beta=\pa_{\o}(\pa_{\o}u) &\mbox{and}& \pa_tu:= \a-\frac12 \tr(\pa^2_{\o\o}u).
\eea
\end{defn}
\no We remark that the path derivatives, whenever they exist, are unique in ``$\dbP_0$-a.s." sense.

\begin{rem}[{\bf Functional It\^o formula}]
\label{rem-Ito1}
{\rm 
When the path derivatives $\pa_t u, \pa_\o u, \pa^2_{\o\o}u$ exist, we have $\a=\pa_tu+\frac12\tr( \pa^2_{\o\o}u)$ and $\beta=\pa_\o u$. In other words, the functional It\^o formula  holds: for
$ t\ge 0$ 
\bea
\label{Ito}
 u_t=u_0+\int_0^t[\pa_t u+\frac12 \tr(\pa^2_{\o\o}u)](s,\cd) ds+\int_0^t \pa_\o u(s,\cd)\cd dB_s,~~\mbox{$\hP_0$-a.s.}
 \eea
 In particular, this implies that $u$ is continuous in $t$. Equivalently, since $\beta=\pa_\o u$ is a semi-martingale, by using the Stratonovich integral, denoted by $\circ dB_s$, one
 has
 \bea
\label{Ito2}
 u_t=u_0+\int_0^t\pa_t u(s,\cd) ds+\int_0^t \pa_\o u(s,\cd)\circ dB_s,\q t\ge 0,~\mbox{$\hP_0$-a.s.}
 \eea
 \qed}
 \end{rem}
\begin{rem}
\label{rem-Ito2}
{\rm 
 (i) The main result in \cite{Dupire} and \cite{CF} is the functional It\^{o} formula (\ref{Ito}), and in \cite{ETZ} the functional It\^{o} formula (\ref{Ito}) is used to define the derivatives. In this sense, our definition is consistent with theirs. 

(ii) In \cite{Dupire} and \cite{CF}, one needs to extend the processes from $\O$ to the space of {\cad} paths. In \cite{ETZ} the definition
is restricted to the space $\O$ only, but it still requires the processes and all the derivatives involved be continuous in $\o$. Our  path-derivatives do not require such regularity, in particular our derivatives are defined only in ``$\dbP_0$-a.s." sense. In this aspect our definition is weaker, and  is convenient for our study of SPDEs in \cite{BMZ}, as typically one cannot expect the solution of a SPDE to be continuous in $\o$. 

(iii) In \cite{Dupire}, \cite{CF}  and \cite{ETZ},  the path derivative $\pa_\o u$ is not required to be an It\^{o} process. In this sense our definition is stronger. This is mainly because our pathwise Taylor expansion below requires stronger regularity than the functional It\^{o} formula.  

(iv) When $u(t,\o) = v(t,\o_t)$ with $v\in C^{1,2}([0,\infty)\times \dbR^d)$, by the standard It\^{o} formula we see that $\pa_\o u(t,\o_t) = \pa_x v(t,\o_t)$. If we assume further that $\pa_x v \in C^{1,2}([0,\infty)\times \dbR^d)$, then $\pa^2_{\o\o} u(t,\o) = \pa^2_{xx} v(t,\o_t)$ and $\pa_t u(t,\o) = \pa_t v(t,\o_t)$. So our  path derivatives are consistent with the standard derivatives in Markovian case. However, as pointed out in (iii), we need a slightly stronger regularity requirements. 
\qed}
\end{rem}

\begin{rem}
\label{rem-noncommute}
{\rm (i) In general the differential operators $\pa_t$ and  $\pa_\o$ cannot commute. Moreover, in the case $d>1$, the Hessian matrix $\pa^2_{\o\o} u$ may not be symmetric,   which implies that in general $\pa_{\o^i}$ and $\pa_{\o^j}$ do not commute either. See Example \ref{eg-noncommute} below. 

(ii) When $u(t,\o) = v(t, \o_{t_1\wedge t},\cds, \o_{t_n\wedge t})$ for some $t_1,\cds, t_n$ and  some deterministic smooth function $v$. Then $\pa^2_{\o\o} u$ is symmetric and  $\pa_{\o^i}$ and $\pa_{\o^j}$ commute. 

(iii) Under the conditions of \cite{Dupire} (or \cite{CF})  
%the notation $\pa^2_{\o\o} u$ is by definition 
the ``Hessian"  $\pa^2_{\o\o} u$ is always symmetric. In fact, being uniformly continuous in $(t,\o)$, the process $u$ in \cite{CF} and
  \cite{Dupire} can be approximated by processes in the form of (ii) above.

(iv) In  \cite{ETZ}  the  $\pa^2_{\o\o} u$ is by definition symmetric. Indeed, the  $\pa^2_{\o\o} u$ in \cite{ETZ}   corresponds to ${1\over 2}[\pa^2_{\o\o} u + (\pa^2_{\o\o} u)^T]$ here. 
However, in this case the relation 
%such a definition does not verify the fact 
$\pa^2_{\o\o}u = \pa_\o(\pa_\o u)$ may fail to hold, which not only  is somewhat unnatural, but also makes the
definition of higher order derivatives much more difficult.  Our new definition modified this point.   
We should also note that the process $u$ in  
Example \ref{eg-noncommute} (ii) below is not continuous in $(t,\o)$, and thus  is not in the framework of \cite{CF},  \cite{Dupire},  or \cite{ETZ}.}
\qed

\end{rem}
\begin{eg}
\label{eg-noncommute}
(i) Let $d=1$ and $ d u = B_t dt$. Then $\pa_\o u = 0$, $\pa_t u = B_t$. It is clear that $\pa_t \pa_\o u = 0 \neq 1= \pa_\o\pa_t u$.

(ii)  Let $d = 2$ and $d u = B^2_t dB^1_t$. Then  $\pa_\o u = [B^2, 0]^T$, and thus $\pa^2_{\o\o} u = \Big[\ba{lll} 0~~1\\0~~0\ea\Big]$ is not symmetric. In particular, $\pa_{\o^1}\pa_{\o^2} u = 0 \neq 1 = \pa_{\o^2}\pa_{\o^1} u$.
\end{eg}

We note that the path derivatives can be extended to any order in a natural way. We now introduce some $L^p$ spaces that will 
be frequently used in the paper. We begin by introducing the following norms on $u\in \dbL^0(\L)$:
% first define $\|u\|_{0,p, T} := \|u\|_{p, T}:=\hE\big[\sup_{0\le t\le T}|u_t|^p\big]^{1/p}$. Then we define recursively:
\bea
\label{normnpT}
\|u\|_{0,p, T} &:=& \|u\|_{p, T}\; :=\;\sup_{0\le t\le T}\Big(\hE\big[|u_t|^p\big]\Big)^{1/p}\\
 \|u\|_{1,p, T} &:=& \|u\|_{p, T} + \sum_{i=1}^{d} \|\pa_{\o^i}u\|_{0,p,T} ;\\
\|u\|_{n,p, T} &:= &\|u\|_{p, T} +\|\pa_t u\|_{n-2, p, T} +  \sum_{i=1}^{d} \|\pa_{\o^i}u\|_{n-1,p,T},\q n\ge 2. \nonumber
\eea
We then define the spaces:  %$ \cH^{[0]}_p(\L) = \dbL^p(\L) $,  and
% \bea
% \label{normKp}
 %\eea 
%Now for $n\ge 1$ we define 
\bea
 \label{cHn}
& \cH^{[n]}_p(\L) := \Big\{u\in \dbL^0(\L): \|u\|_{n,p, T} <\infty,  \forall T > 0\Big\}, ~~n\ge 0; 
\eea
%The space $\cH^{[n]}_p(\L)$ can be extended to the  space. 
We shall also define the following H\"older norms: for any $\a \in (0,1)$,
 \bea
 \label{cHn2}
\|u\|_{0,p,\a, T} &:= & \|u\|_{p,\a,T}\;:=\;  \|u\|_{p, T} + \hE\Big[\sup_{0\le t \le T, \d>0}{|u_{t+\d}  - u_t|^p \over \d^{p\a\over 2}}\Big]^{1/p}, \nonumber\\
 \|u\|_{1,p, \a, T} &:=& \|u\|_{p,\a, T} +  \sum_{i=1}^{d} \|\pa_{\o^i}u\|_{p,\a, T},\\
\|u\|_{n,p, \a, T} &:= &\|u\|_{p, \a, T} +\|\pa_t u\|_{n-2, p, \a, T} + \sum_{i=1}^{d} \|\pa_{\o^i}u\|_{n-1,p,\a, T}, \q n\ge 2. \nonumber
\eea
Then  we define correspondingly:
\bea
\cH^{[n]+\a}_p(\L) := \Big\{u\in \cH^{[n]}(\L): \|u\|_{n,p, \a,T} <\infty,  \forall T > 0\Big\}, ~~n\ge 0.
 \eea 

 It should be noted that $\cH^{[n+1]}_p(\L)$  is not a subspace of  $\cH^{[n]+\a}_p(\L)$(!). However, one can show that
if $p>0$ is large enough (more precisely $p> {2\over 1-\a}$), then it holds that  $\cH^{[n+2]}_p(\L) \subset \cH^{[n]+\a}_p(\L)$ for any 
$n\ge 0$ (see Lemma \ref{lem-inclusion} below).
%Clearly the path derivatives can be extended to any order. We then define 
% \bea
% \label{normKp}
%\dbL^p(\L) := \Big\{u\in \dbL^0(\L):   \|u\|_{p,T} < \infty , \forall T > 0\Big\} &\mbox{where}& \|u\|_{p,T}^p:=\sup_{0\le t\le T}\hE[|u_t|^p],
% \eea 
% and
% \bea
% \label{cHn}
%& \cH^{[n]}_p(\L) := \Big\{u\in \dbL^0(\L): \|u\|_{n,p, T} <\infty,  \forall T > 0\Big\}, ~~n\ge 0; \q \mbox{where}&\\
%&\|u\|_{0,p, T} := \|u\|_{p, T},\q \|u\|_{1,p, N} := \|u\|_{p, N} + \sum_{i=1}^{d} \|\pa_{\o^i}u\|_{0,p,T},& \nonumber\\
%&\|u\|_{n,p, N} := \|u\|_{p, N} +\|\pa_t u\|_{n-2, p, T} +  \sum_{i=1}^{d} \|\pa_{\o^i}u\|_{n-1,p,T},\q n\ge 2.& \nonumber
%\eea

\subsection{Path derivatives for random fields}
  Let  $\cO\subset \hR^{d'}$ be an %convex 
  open domain,  $Q:= [0,\infty)\times \cO$, and $\hat \L := Q \times \O$. We denote by $\dbL^0(\hat\L, \dbR^{m\times n})$ the set of  $\dbF$-progressively  measurable random fields $ u : \hat\L \to \dbR^{m\times n}$, and $\dbL^0(\hat\L):= \dbL^0(\hat\L, \dbR)$.  When there is no confusion, we shall omit the variable $\o$ in $u$ and  write it as $u(t,x)$.  Given $u\in \dbL^0(\hat\L)$, we define its derivatives in $x$ in the standard way, and for any fixed $x$, its  path-derivatives in $(t,\o)$ in the sense of Definition \ref{Dw}.

Notice that $Q$ is not compact.   For any $N\ge 1$, denote
\bea
\label{KN}
K_N := \Big\{x\in \cO: |x|\le N, ~d(x, \cO^c) \ge {1\over N}\Big\} &\mbox{and}& Q_N := [0, N] \times K_N.
\eea
It is clear that $K_N$ is %convex, 
compact, increasing in $N$, and
\bea
\label{KN2}
\bigcup_{N=1}^\infty K_N = \cO, ~\bigcup_{N=1}^\infty Q_N = Q,~\mbox{and}~ x+h\in K_{N+1}~\mbox{for any}~x\in K_N, |h|\le {1\over N(N+1)}.
\eea
 Similar to the process case we can define the norms:  for $p\ge 1$ and $n\ge 2$,
%We then define 
 \bea
 \label{normKp-space}
%\left\{\ba{lllll}
 \|u\|_{0,p, N}&:=& \|u\|_{p, N}:=\sup_{t\in [0, N]}\hE\Big[\sup_{x\in K_N}|u(t,x)|^p\Big]^{1/p}, 
 % + \sup_{x\in K_N}\hE\Big[\sup_{t\in [0,N]}|u(t,x)|^p\Big]^{1/p},
\nonumber\\
\|u\|_{1,p, N}& :=& \|u\|_{p, N} + \sum_{i=1}^{d'} \|\pa_{x_i} u\|_{0,p,N} + \sum_{i=1}^{d} \|\pa_{\o^i}u\|_{0,p,N}, \nonumber\\
 \|u\|_{n,p, N} &:=& \|u\|_{p, N} +\|\pa_t u\|_{n-2, p, N} + \sum_{i=1}^{d'} \|\pa_{x_i} u\|_{n-1,p,N} + \sum_{i=1}^{d} \|\pa_{\o^i}u\|_{n-1,p,N},
 % ~n\ge 2. 
\nonumber \\
% \label{normKp2-space}
\|u\|_{0,p,\a, N} &:= &\|u\|_{p, \a, N}:=  \|u\|_{p, N} + \sup_{x\in K_N}\hE\big[\sup_{0\le t < t'\le N}{|u(t,x) - u(t',x)|^p\over |t-t'|^{{p\a\over 2}} }\big]^{1/p}\\
 &&+\sup_{0\le t\le N}\hE\Big[\sup_{x, x' \in K_N}{|u(t,x) - u(t,x')|^p\over|x-x'|^{p\a}} \Big]^{1/p}, \qq \a\in(0,1),\nonumber\\
 \|u\|_{1,p, \a, N}& := &\|u\|_{p,\a, N} + \sum_{i=1}^{d'} \|\pa_{x_i} u\|_{p,\a, N} + \sum_{i=1}^{d} \|\pa_{\o^i}u\|_{p,\a, N}, \nonumber\\
\|u\|_{n,p, \a, N} &:=& \|u\|_{p, \a, N} +\|\pa_t u\|_{n-2, p, \a, N} + \sum_{i=1}^{d'} \|\pa_{x_i} u\|_{n-1,p,\a, N} + \sum_{i=1}^{d} \|\pa_{\o^i}u\|_{n-1,p,\a, N}. \nonumber
\eea 
We  now define the following spaces: for $p\ge 1$, $0<\a<1$, and $n\ge 0$, 
\bea
 \label{cHn-space}
 \left.\ba{lll}
%\cH^{[0]}_p(\hat\L) &:=&\dbL^p(\hat\L) := \Big\{u\in \dbL^0(\hat\L):   \|u\|_{p,N} < \infty , \forall N > 0\Big\},\nonumber\\
\dis \cH^{[n]}_p(\hat\L) := \Big\{u\in \dbL^0(\hat\L): \|u\|_{n,p, N} <\infty,  \forall N > 0\Big\}, \ms \\
% \label{cHn2-space}
\dis \cH^{[n]+\a}_p (\hat\L):= \Big\{u\in \cH^{[n]}_p(\hat\L): \|u\|_{n,p, \a,N} <\infty,  \forall N > 0\Big\}.
\ea\right.
\eea
%for $n\ge 2$ in the last line.  
Again, as we shall see in Lemma \ref{lem-inclusion}, one can show that
\bea
\label{inclusion}
\mbox{$\cH^{[n+2]}_p(\hat\L) \subset \cH^{[n]+\a}_p(\hat\L)$ for any $n\ge 0$,  $\a\in (0,1)$ and $p> {2\over 1-\a}$.}
\eea

\subsection{It\^o random fields and It\^{o}-Ventzell Formula}
\ms
  We recall that $u \in \dbL^0(\hat\L)$ is called an It\^o random field if, for any $x\in \cO$,
\bea
\label{uito2}
u(t,x) = u_0(x) + \int_0^t \a(s,x) ds+\int_0^t \b(s,x) \cd dB_s, \q  t\ge 0, ~\hP_0\mbox{-a.s.}
\eea
where $\a \in \dbL^0(\hat\L)$, $\b \in  \dbL^0(\hat\L, \dbR^d))$ satisfy $\int_0^t [|\a(s,x)|+|\b(s,x)|^2]ds <\infty$, $\dbP_0$-a.s. for all $(t,x)\in Q$. 

It is worth noting that, in contrast to Remark \ref{rem-noncommute},  the spatial derivative $\pa_x$ commutes with 
both $\pa_t$ and $\pa_\o$. In fact, we have the following result. Since the proof is quite straightforward, we omit it. 
\begin{lem}
\label{lem-commute} Let $u \in \dbL^0(\hat\L)$ be an It\^o random field in the form of (\ref{uito2}).

(i) %Assume $\pa_x u_0$, $\pa_x \a \in \dbL^1(\L, \dbR^{d'})$, $\pa_x \b\in \dbL^2(\L, \dbR^{d\times d'})$ exist, 
Assume $u_0, \a, \b$ are differentiable in $x$, and for any $N>0$, the processes $|\pa_x \a(\cd, x,\cd)|$ and $|\pa_x \b(\cd, x,\cd)|^2$ are uniformly integrable on $[0, N]\times \O$, uniformly on $x\in K_N$. Then $\pa_x u$ exists and is also an It\^{o} random field: for each $i=1,\cds, d'$,
\beaa
\pa_{x_i} u(t,x) =\pa_{x_i} u_0(x) + \int_0^t \pa_{x_i} \a(s,x) ds+\int_0^t \pa_{x_i} \b(s,x) \cd dB_s, \q  t\ge 0, ~\hP_0\mbox{-a.s.}
\eeaa
In particular, this implies that
\bea
\label{commute1}
\pa_\o \pa_x u = \pa_x \pa_\o u.
\eea
(ii)  Assume further that $\b$ is an It\^{o} random field and each of its components satisfies the property of $u$ in (i), then 
\bea
\label{commute2}
\pa^2_{\o\o} \pa_x u = \pa_x \pa^2_{\o\o} u &\mbox{and}& \pa_t \pa_x u = \pa_x \pa_t u.
\eea
\end{lem}

As an important application of the path derivatives,  we recast the It\^o-Ventzell formula, which turns out to be exactly the same as a multidimensional funtional It\^{o} formula.
 \begin{prop}[{\bf It\^o-Ventzell formula}]
 \label{prop-ItoVentzell}
 Let $X \in \cH^{[2]}_2(\L)$ taking values in $\cO$ and $u\in \cH^{[2]}_2(\hat\L)$  such that $u(\cd,\o) \in C^{0,2}(Q)$ and $\pa_\o u(\cd,\o) \in C^{0,1}(Q)$, for $\dbP_0$-a.e. $\o$.
%$$d X_t = b_t dt + \si_t dB_t, \qq t\ge 0.$$
Then %for any $u\in \cH^{[2]}_2(\L)$, 
the following chain rule for our path derivatives holds:
 \bea
 \label{chainrule}
 \pa_t [u(t, X_t, \o)]  = \pa_t u + \pa_x u \cd \pa_t X_t;  &&  \pa_{\o} [u(t, X_t, \o)]  =  \pa_\o u +  [\pa_\o X_t]^T \pa_x u.
  \eea
In particular, if  
$d X_t = b_t dt + \si_t \cd dB_t$, $d u(t,x,\o) = \a(t,x) dt + \b(t,x) \cd dB_t$, $t\ge 0$, 
then the It\^o-Ventzell formula holds:
%following identity (suppressing variables) holds: 
for $t\ge 0$, $\hP_0$-a.s.,
 \bea
 \label{Ito-Ventzell}
 d u(t, X_t, \o) &=& \big[\pa_t u + \pa_x u \cd \pa_t X_t  \big]dt + \big[ \pa_\o u +  [\pa_\o X_t]^T \pa_x u\big] \circ dB_t\\
 &=& \big[ \a  + \pa_x u \cd b_t + {1\over 2} \pa_{xx}^2 u : \si_t\si^T_t +\pa_x \b : \si_t\big] dt + [\b+\si_t^T \pa_x u] \cd dB_t.
 \nonumber
 \eea
 \end{prop}

\no{\it Proof.} Since $X \in \cH^{[2]}_2(\L)$ and $u\in \cH^{[2]}_2(\hat\L)$, 
%all the (path) derivatives exist, and 
one can write
$$d X_t = b_t dt + \si_t \cd dB_t, \q\mbox{ and } \q d u(t,x,\o) = \a(t,x) dt + \b(t,x) \cd dB_t,$$
 where $b^i = \pa_t X^i + {1\over 2} \tr(\pa^2_{\o\o} X^i)$, $\si^i = \pa_\o X^i$, $i=1,\cds d'$, and $\a = \pa_t u + {1\over 2}\tr(\pa^2_{\o\o} u)$, $\b = \pa_\o u$. 

Next,   under our conditions we may apply the   standard It\^{o}-Ventzell formula and obtain
%from (\ref{cHn-space}) we see that $u\in\cH^{[2]}_2(\hat \L)$ implies that $u(\cd,\cd, \o)\in C^{1,2}([0,T]\times\cO)$ for $\hP$-a.s., $\o$, applying the standard It\^{o}-Ventzell formula we have
\beaa
&d u(t, X_t) = \big[ \a  + \pa_x u \cd b_t + {1\over 2} \pa_{xx}^2 u : \si_t\si^T_t +\pa_x \b : \si_t\big](t, X_t) dt + [\b+\si_t^T \pa_x u] (t, X_t)\cd dB_t.&
\eeaa
Therefore the definition of path derivatives leads to that
\bea
\label{chainrule0}
\left\{\ba{lll}
\pa_{\o} [u(t, X_t)]  =  \b+\si_t^T \pa_x u =  (\pa_\o u)(t,X_t) +  [\pa_\o X_t]^T[ (\pa_x u)(t, X_t) ];\ms\\
\dis \pa_t [u(t, X_t, \o)] =  \big[ \a  + \pa_x u \cd b_t + {1\over 2} \pa_{xx}^2 u : \si_t\si^T_t +\pa_x \b : \si_t\big]  - 
 {1\over 2} \tr\big(\pa_{\o\o}[u(t, X_t)]\big).
 \ea\right.
 \eea
Note that $\pa^2_{\o\o}u =\pa_\o[\pa_\o u]$, differentiating $\pa_\o u(t, X_t)$ again we have
\beaa
\pa^2_{\o\o}[u(t, X_t)] = \pa_{\o\o}u  + \pa_\o X_t [\pa_{x\o} u]^T +   \sum_{i=1}^{d'} \big[\pa_{\o\o} X^i_t \pa_{x_i} u + \pa_\o X^i_t [\pa_{x_i \o} u + \pa_{x_i x}u \pa_\o X^i]^T\big].
\eeaa
Now plugging this into (\ref{chainrule0}) and recalling the definition of $\pa_\o X$, $\pa^2_{\o\o}X$, $\pa_\o u$, $\pa^2_{\o\o} u$, 
with some simple computation we prove 
(\ref{chainrule}), whence (\ref{Ito-Ventzell}),  immediately.
  \qed

\begin{rem}
\label{ChainRule}
{\rm (i) If $u$ is deterministic, then $\beta=\pa_\o u=0$, and we have the It\^o formula.

\no (ii) As the ``chain rule" (\ref{chainrule})  completely characterizes the expression (\ref{Ito-Ventzell}), we may refer to it as ``pathwise It\^o-Ventzell formula".}
\qed
\end{rem}

% as  determines the 
\subsection{Multiple differentiation and integration}
\label{sect-notation}
Our Taylor expansion will involve multiple differentiation and integration.  However, due to the noncommutative property of the path derivatives in Remark \ref{rem-noncommute} and Example \ref{eg-noncommute}, we need to specify the differentiation and integration indices precisely. To simplify presentation, we first introduce some notations. For $i = 0,1,\cds, d$, define
\bea
\label{I}
\left.\ba{lll}
\pa_i u :=\pa_t u,\q u_t d_i t :=  u_t dt,\q \qq &\mbox{if}~i=0;\\
 \pa_i u := \pa_{\o^i} u,~~ u_t d_i t := u_t \circ dB^i_t, &\mbox{if}~1\le i\le d.
\ea\right.
\eea
Next, for $\th = (\th_1,\tilde \th) =(\th_1, \th_2,\cds, \th_n) \in \{ 0,1,\cds, d\}^n$ and $s<t$, %in light of (\ref{Taylor+}) and (\ref{R+}), 
we define recursively by:
\bea
\label{DI}
 \cD^\th_\o u := \pa_{\th_1}(\cD_\o^{\tilde\th} u),\q \cI^\th_{s, t}(u) := \cI^{\tilde\th}_{s, t}\Big(\int_s^\cd u_r d_{\th_1}r\Big),\q \cI^\th_{s, t} := \cI^\th_{s, t}(1) .
\eea
Notice that the above definition also implies, for $\th = (\bar\th, \th_n)$,
\bea
\label{DI2}
\left.\ba{c}
\dis \cD^\th_\o u = \pa_{\th_1} \cds \pa_{\th_n}  u = \cD_\o^{\tilde\th}(\pa_{\th_n} u),\\
\dis  \cI^\th_{s, t}(u) :=  \int_s^t \int_s^{t_{n}}\cds\int_s^{t_2} u_{t_1} d_{\th_1} t_1  \cds d_{\th_n}  t_n = \int_s^t \cI^{\bar\th}_{s, r}(u) d_{\th_n}r.
\ea\right.
\eea
Moreover, for the purpose of backward expansion later, we introduce
%notice that the terms in (\ref{Taylor-}) and (\ref{R-}) reverse the order of integration of the corresponding terms in (\ref{Taylor+}) and (\ref{R+}), we define:
\bea
\label{-th}
-\th := (\th_n,\cds,\th_1) &\mbox{and}& \cI^\th_{t, s}(u) := (-1)^n \cI^{-\th}_{s, t}(u) ~\mbox{for}~s<t.
\eea

Noting the relation between the horizontal derivative $\pa_t u$ and  $\pa^2_{\o\o} u$ (cf. (\ref{Dtu})), we introduce the
following ``weighted norm":  for $\th \in \{0,1,\cds, d\}^n$,
\bea
\label{lnorm}
|\th|_0 := n,\q |\th| := n + \sum_{i=1}^n 1_{\{\th_i = 0\}}.
\eea
Moreover, when $|\th|=0$,  we take the notational convention that
\bea
\label{th=0}
\cD^\th_\o u := u,\q \cI^\th_{s,t}(u) := u_{t}.
\eea
Due to the commutative property of Lemma \ref{lem-commute}, the high order differentiation operator in $x$ is simpler.  
Let $\dbN$ be the set of nonnegative integers. For $\ell = (\ell_1,\cds, \ell_{d'}) \in \dbN^{d'}$, denote:
\bea
\label{Dx}
\cD^\ell_x u := \pa^{\ell_1}_{x_1}\cds\pa^{\ell_{d'}}_{x_{d'}}u,\q  x^\ell := \Pi_{i=1}^{d'} x_i^{\ell_i},\q \ell ! := \Pi_{i=1}^{d'} \ell_i!,\q |\ell| := \sum_{i=1}^{d'} \ell_i
\eea
%Similar to (\ref{th=0}), we 
We shall %denote  $|\ell| := \sum_{i=1}^{d'} \ell_i$, and 
set $\cD^\ell_x u := u$, $x^\ell := 1$, and $\ell! := 1$, if $|\ell|=0$. 

Furthermore, together with (\ref{lnorm}), we can  introduce a ``weighted norm" on the index set
$\Th:= \cup_{n=0}^\infty \{0,1,\cds, d\}^n \times \dbN^{d'}
$:
\bea
\label{lnorm-space}
 |(\th, \ell)| := |\th|+|\ell|, ~ \q \forall (\th,\ell) \in \Th.
\eea
%We now consider the following  
Note that if we denote $ \Th_n := \{(\th, \ell)\in \Th: |(\th,\ell)|\le n\}$, then by applying Lemma \ref{lem-commute} one can easily check
 that: 
%\bea
%\label{cHn+1}
%&\mbox{
if  $u\in \cH^{[n]}_2$, then 
%}&\nonumber\\
%&\!\!\!\!\mbox{
all derivatives of $u$ up to order $n$ can be written as $\cD^\ell_x \cD^\th_\o u$ for some $(\th, \ell)\in \Th_n$ (counting ``$\pa_t$" as 
a second order derivative!).

%\eea

\section{Taylor Expansion for It\^{o} Processes (Second Order Case)}
\setcounter{equation}{0}

%\section{Heuristic derivation of second order  expansion} 
%Due to the non-commutative properties of our path derivatives, we need to introduce some notations for multiple differentiation and integration in high dimensions. To motivate them, 

In this section we give some heuristic arguments for
the simplest second order Taylor expansion for It\^o processes. 
%We begin with 
%of a process $u$ which is smooth enough in our sense. 
We shall establish both forward and backward temporal expansions.
% and the same idea will be applied to the general cases in the next section. 
%We should note that the  backward expansion is 
%particularly important in the study of viscosity solution of SPDEs 
%%will rely on the backward temporal expansion, 
%(see our accompanying paper \cite{BMZ}). 

In what follows we shall always denote, for $s<t$,
%\bea
%\label{fst}
$\f_{s,t} := \f_t - \f_s$,
%\eea
and we will use the following simple fact frequently:   for any  semimartingales $\xi, \eta, \g$, 
\bea
\label{composition}
\int \xi_t \circ (\eta_t \circ d\g_t) = \int (\xi_t \eta_t) \circ d\g_t =\int \eta_t \circ (\xi_t \circ d\g_t).
\eea

\subsection{Forward Temporal Expansion.}
 Let $t\ge 0$, $\d>0$, and denote $t_\d := t+\d$. Repeatedly applying the functional It\^{o} formula (\ref{Ito2}) formally we have
\beaa
u_{t_\d}  &=& u_t + \int_t^{t_\d} \pa_t u_s ds + \sum_{i=1}^d \int_t^{t_\d} \pa_{\o^i} u_s \circ dB^i_s\\
&=& u_t + \sum_{i=1}^d \pa_{\o^i} u_t B^{i}_{t, t_\d} +  \int_t^{t_\d} \pa_t u_s ds + \sum_{i=1}^d \int_t^{t_\d} [\pa_{\o^i} u]_{t,s} \circ dB^i_s \\
&=& u_t + \sum_{i=1}^d \pa_{\o^i} u_t B^{i}_{t,t_\d} +  \int_t^{t_\d} \pa_t u_s ds \\
&&\!\!\!+ \sum_{i=1}^d \int_t^{t_\d}\Big[\int_t^s \pa_{t\o^i} u_r dr +  \sum_{j=1}^d \int_t^s \pa_{\o^j\o^i} u_r \circ dB^j_r\Big]\circ dB^i_s\\
&=&u_t + \sum_{i=1}^d \pa_{\o^i} u_t B^{i}_{t,t_\d} +  \pa_t u_t \d + \sum_{i, j=1}^d   \pa_{\o^j\o^i} u_t \int_t^{t_\d} B^{j}_{t,s} \circ dB^i_s  \\
&&\!\!\!+ \int_t^{t_\d} [\pa_t u]_{t,s} ds + \sum_{i=1}^d \int_t^{t_\d}\int_t^s \pa_{t\o^i} u_r dr \circ dB^i_s +  \sum_{i,j=1}^d \int_t^{t_\d}\int_t^s [\pa_{\o^j\o^i} u]_{t,r} \circ dB^j_r\Big]\circ dB^i_s.
\eeaa
Here we used the fact that $\pa_{\o^i} u_t$ and $\pa_{\o^j\o^i}u_t$ are $\cF_t$-measurable and can be moved out from the related stochastic integrals. (We note that this will not be the case when we consider backward temporal expansion later.) Then
\bea
\label{Taylor+}
u_{t_\d}&=& u_t + \sum_{i=1}^d \pa_{\o^i} u_t B^{i}_{t,t_\d} +  \pa_t u_t \d + \sum_{i, j=1}^d   \pa_{\o^j\o^i} u_t \int_t^{t_\d} B^{j}_{t,s} \circ dB^i_s  + R_2(t,\d),
\eea  
where
\bea
\label{R+}
 R_2(t,\d) &:=& \int_t^{t_\d} \int_t^s \pa_{tt} u_rdrds + \sum_{i,j=1}^d \int_t^{t_\d} \int_t^s \int_t^r \pa_{t\o^j\o^i} u_\k d\k \circ dB^j_r \circ dB^i_s\\
 &&+ \sum_{i=1}^d \int_t^{t_\d} \int_t^s\pa_{\o^it} u_r \circ dB^i_r ds + \sum_{i=1}^d \int_t^{t_\d} \int_t^s\pa_{t\o^i} u_r dr \circ dB^i_s\nonumber\\
 && +  \sum_{i,j,k=1}^d \int_t^{t_\d}\int_t^s \int_t^r \pa_{\o^k\o^j\o^i} u_\k \circ dB^k_\k \circ dB^j_r\circ dB^i_s. \nonumber
\eea

To simplify the presentations let us make use of the notations for multiple derivatives and integrations defined in 
(\ref{I})--(\ref{DI2}). 
Then it is straightforward to check that (\ref{Taylor+}) and (\ref{R+}) can be rewritten as a more compact form:
\bea
\label{Taylor+2}
u_{t_\d}&=&  u_t + \sum_{i=0}^d \cD^{(i)}_\o u_t \cI^{(i)}_{t,t_\d}  + \sum_{i, j=1}^d   \cD^{(j,i)}_\o u_t \cI^{(j,i)}_{t,t_\d}   + R_2(t,\d)\\
\label{R+2}
 R_2(t,\d) &:=& \cI^{(0,0)}_{t,t_\d}(\cD^{(0,0)}_\o u) + \sum_{i,j=1}^d  \cI^{(0,j,i)}_{t,t_\d}(\cD^{(0, j,i)}_\o u)\\
 &&+ \sum_{i=1}^d  \cI^{(i,0)}_{t,t_\d}(\cD^{(i,0)}_\o u)  + \sum_{i=1}^d \cI^{(0,i)}_{t,t_\d}(\cD^{(0,i)}_\o u) +  \sum_{i,j,k=1}^d \cI^{(k,j,i)}_{t,t_\d}(\cD^{(k,j,i)}_\o u). \nonumber
 \eea

\subsection{Backward Temporal Expansion.}
 Let $0<\d\le t$, and denote $t^-_\d := t-\d$. Then similar to the forward expansion we can obtain 
\bea
\label{ut-d}
u_{t^-_\d} &=& u_t - \int_{t^-_\d}^t \pa_t u_s ds - \sum_{i=1}^d\int_{t^-_\d}^t \pa_\o u^i_s \circ dB^i_s\\
&=&  u_t - \pa_t u_{t^-_\d} \d  - \int_{t^-_\d}^t [\pa_t u]_{t^-_\d,s} ds - \sum_{i=1}^d\Big[\pa_{\o^i} u_{t^-_\d} B^i_{t^-_\d,t} 
+ \int_{t^-_\d}^t [\pa_{\o^i} u]_{t^-_\d, s} \circ dB^i_s\Big]. \nonumber
%\\
%&=& u_t -  \sum_{i=1}^d \pa_{\o^i} u_t B^i_{t^-_\d,t} +   \sum_{i=1}^d[\pa_{\o^i} u]_{t^-_\d, t} B^i_{t^-_\d,t} - \int_{t^-_\d}^t \pa_t u_s ds -  \sum_{i=1}^d\int_{t^-_\d}^t [\pa_{\o^i} u]_{t^-_\d, s} \circ dB^i_s.\\
%&&
%+  \sum_{i,j=1}^d  \pa_{\o^j\o^i} u_{t^-_\d} \int_{t^-_\d}^tB^i_{t^-_\d, s}  \circ dB^j_s\\
\eea
We should note that  the above  expansion would be around $t^-_\d$ instead of $t$, we therefore modify it as follows. 
First, we write
\bea
\label{t-d-to-t}
\pa_t u_{t^-_\d} \d=\pa_t u_t\d-[\pa_t u]_{t^-_\d,t} \d,\q
 \pa_{\o^i} u_{t^-_\d}B^i_{t^-_\d,t}=\pa_{\o^i} u_t B^i_{t^-_\d,t}-[\pa_{\o^i} u]_{t^-_\d,t} B^i_{t^-_\d,t}.
 \eea
 Next, we apply integration by parts formula and/or (standard) It\^{o} formula to get
\bea
\label{IBP1}
 &&[\pa_t u]_{t^-_\d, t} \d - \int_{t^-_\d}^t [\pa_t u]_{t^-_\d, s} ds =\int_{t^-_\d}^ t (s-t^-_\d) d(\pa_t u_s)\nonumber\\
 &=& \int_{t^-_\d}^ t \pa_{tt} u_s (s-t^-_\d) ds +  \sum_{i=1}^d \int_{t^-_\d}^ t \pa_{\o^i t} u_s (s-t^-_\d) \circ dB^i_s;\nonumber\\
&&[\pa_{\o^i} u]_{t^-_\d, t} B^i_{t^-_\d,t} - \int_{t^-_\d}^t [\pa_{\o^i} u]_{t^-_\d, s} \circ dB^i_s \\
&=& \int_{t^-_\d}^t B^i_{t^-_\d, s} \circ  d(\pa_{\o^i} u_s ) = \int_{t^-_\d}^t B^i_{t^-_\d, s}  \circ \Big(\pa_{t\o^i} u_s ds +\sum_{j=1}^d \pa_{\o^j\o^i} u_s \circ dB^j_s\Big) \nonumber\\
&=&  \int_{t^-_\d}^t \pa_{t\o^i} u_s B^i_{t^-_\d, s}  ds + \sum_{j=1}^d  \int_{t^-_\d}^t \pa_{\o^j\o^i} u_s B^i_{t^-_\d, s}  \circ dB^j_s\nonumber\\
&=& \int_{t^-_\d}^t \pa_{t\o^i} u_s B^i_{t^-_\d, s}  ds +\sum_{j=1}^d \pa_{\o^j\o^i} u_t  \int_{t^-_\d}^t  B^i_{t^-_\d, s}  \circ dB^j_s\nonumber\\
&&-\sum_{j=1}^d  [\pa_{\o^j\o^i}]_{t_\d^-, t}  \int_{t^-_\d}^t  B^i_{t^-_\d, s}  \circ dB^j_s+ \sum_{j=1}^d  \int_{t^-_\d}^t [\pa_{\o^j\o^i} u]_{t_\d^-,s} B^i_{t^-_\d, s}  \circ dB^j_s,\nonumber
 \eea
and 
%   our purpose is to find the expansion around the time $t$, thus the Apply integration by parts formula,
   \bea
   \label{IBP2}
&&[\pa_{\o^j\o^i} u]_{t^-_\d, t} \int_{t^-_\d}^t B^i_{t^-_\d, s}  \circ dB^j_s - \int_{t^-_\d}^t [\pa_{\o^j\o^i} u]_{t^-_\d,s} B^i_{t^-_\d, s}  \circ dB^j_s\\
 &=&  \int_{t^-_\d}^t \big(\int_{t^-_\d}^s B^i_{t^-_\d, r}  \circ dB^j_r\big) \circ d ( \pa_{\o^j\o^i} u_s) \nonumber\\
 &=& \int_{t^-_\d}^t \pa_{t\o^j\o^i} u_s \big(\int_{t^-_\d}^s B^i_{t^-_\d, r}  \circ dB^j_r\big) ds +  \sum_{k=1}^d\int_{t^-_\d}^t \pa_{\o^k \o^j\o^i} u_s \big(\int_{t^-_\d}^s B^i_{t^-_\d, r}  \circ dB^j_r\big) \circ dB^k_s.\nonumber
 \eea
Plugging  (\ref{IBP2}) into \reff{IBP1} and then plugging (\ref{t-d-to-t})--(\ref{IBP1}) into (\ref{ut-d}) we obtain
\bea
\label{Taylor-}
u_{t^-_\d} 
%&=& u_t -  \sum_{i=1}^d \pa_{\o^i} u_t B^i_{t^-_\d,t}  - \int_{t^-_\d}^t \pa_t u_s ds\nonumber\\
%&&+\sum_{i=1}^d  \int_{t^-_\d}^t \pa_{t\o^i} u_s B^i_{t^-_\d, s}  ds + \sum_{i,j=1}^d  \int_{t^-_\d}^t \pa_{\o^j\o^i} u_s B^i_{t^-_\d, s}  \circ dB^j_s\nonumber\\
%&=&  u_t -  \sum_{i=1}^d \pa_{\o^i} u_t B^i_{t^-_\d,t}  - \pa_t u_{t^-_\d} \d +  \sum_{i,j=1}^d  \pa_{\o^j\o^i} u_{t^-_\d} \int_{t^-_\d}^tB^i_{t^-_\d, s}  \circ dB^j_s\\
%&&- \int_{t^-_\d}^t [\pa_t u]_{t^-_\d, s} ds+\sum_{i=1}^d  \int_{t^-_\d}^t \pa_{t\o^i} u_s B^i_{t^-_\d, s}  ds + \sum_{i,j=1}^d  \int_{t^-_\d}^t [\pa_{\o^j\o^i} u]_{t^-_\d,s} B^i_{t^-_\d, s}  \circ dB^j_s\\
%&=&  u_t -  \sum_{i=1}^d \pa_{\o^i} u_t B^i_{t^-_\d,t}  - \pa_t u_{t} \d +  \sum_{i,j=1}^d  \pa_{\o^j\o^i} u_{t} \int_{t^-_\d}^t B^i_{t^-_\d, s}  \circ dB^j_s\\
%&&+ \pa_t u_{t^-_\d, t} \d -  \sum_{i,j=1}^d  \pa_{\o^j\o^i} u_{t^-_\d, t} \int_{t^-_\d}^t B^i_{t^-_\d, s}  \circ dB^j_s\nonumber\\
%&&- \int_{t^-_\d}^t [\pa_t u]_{t^-_\d, s} ds+\sum_{i=1}^d  \int_{t^-_\d}^t \pa_{t\o^i} u_s B^i_{t^-_\d, s}  ds + \sum_{i,j=1}^d  \int_{t^-_\d}^t [\pa_{\o^j\o^i} u]_{t^-_\d,s} B^i_{t^-_\d, s}  \circ dB^j_s\nonumber\\
 &=& u_t -  \sum_{i=1}^d \pa_{\o^i} u_t B^i_{t^-_\d,t}  - \pa_t u_{t} \d +  \sum_{i,j=1}^d  \pa_{\o^j\o^i} u_{t} \int_{t^-_\d}^t B^i_{t^-_\d, s}  \circ dB^j_s + R_2(t, -\d),
 \eea
 where
 \bea
 \label{R-}
 R_2(t, -\d) &=& \int_{t^-_\d}^ t \pa_{tt} u_s (s-t^-_\d) ds  - \sum_{i,j=1}^d \int_{t^-_\d}^t \pa_{t\o^j\o^i} u_s \big(\int_{t^-_\d}^s B^i_{t^-_\d, r}  \circ dB^j_r\big) ds\nonumber\\
 &&+ \sum_{i=1}^d  \int_{t^-_\d}^t \pa_{t\o^i} u_s B^i_{t^-_\d, s}  ds  +  \sum_{i=1}^d \int_{t^-_\d}^ t \pa_{\o^i t} u_s (s-t^-_\d) \circ dB^i_s\\
 && - \sum_{i,j,k=1}^d\int_{t^-_\d}^t \pa_{\o^k \o^j\o^i} u_s \big(\int_{t^-_\d}^s B^i_{t^-_\d, r}  \circ dB^j_r\big)  \circ dB^k_s.\nonumber
 \eea
 Using the notations for multiple derivatives and integrations again, we see that (\ref{Taylor-}) and (\ref{R-})  can again
 be written as the compact form:
\bea
\label{Taylor-2a}
u_{t^-_\d} &=& u_t - \sum_{i=0}^d \cD^{(i)}_\o u_t \cI^{(i)}_{t^-_\d,t} + \sum_{i, j=1}^d   \cD^{(j,i)}_\o u_t \cI^{(i,j)}_{t^-_\d,t}+ R_2(t, -\d);
\eea
where 
\bea
  \label{R-2}
 R_2(t, -\d) &:=& \cI^{(0)}_{t^-_\d, t}\big(\cD^{(0,0)}_\o u_\cd\; \cI^{(0)}_{t^-_\d, \cd}\big) - \sum_{i,j=1}^d \cI^{(0)}_{t^-_\d,t}\big(\cD^{(0,j,i)}_\o u_\cd \; \cI^{(i,j)}_{t^-_\d, \cd}\big)+ \sum_{i=1}^d  \cI^{(0)}_{t^-_\d, t}\big(\cD^{(0,i)}_\o u_\cd\; \cI^{(i)}_{t^-_\d, \cd}\big)\nonumber\\
% &&\\
 &&  + \sum_{i=1}^d  \cI^{(i)}_{t^-_\d, t}\big(\cD^{(i,0)}_\o u_\cd\; \cI^{(0)}_{t^-_\d, \cd}\big)- \sum_{i,j,k=1}^d  \cI^{(k)}_{t^-_\d, t}\big(\cD^{(k,j,i)}_\o u_\cd\; \cI^{(i,j)}_{t^-_\d, \cd}\big).
\eea
We should point out here that (\ref{Taylor-2a}) is slightly different from (\ref{Taylor+2}). But by applying  the 
relation (\ref{-th}) we can rewrite (\ref{Taylor-2a}) as 
\bea
  \label{Taylor-2}
 u_{t^-_\d}&=& u_t + \sum_{i=0}^d \cD^{(i)}_\o u_t \cI^{(i)}_{t,t^-_\d} + \sum_{i, j=1}^d   \cD^{(j,i)}_\o u_t \cI^{(j,i)}_{t,t^-_\d}   + R_2(t,-\d).
 \eea
We see that (\ref{Taylor-2}) is indeed consistent with the forward expansion (\ref{Taylor+2})(!).

\begin{rem}
\label{rem-roughpath}
{\rm (i) If we define, for $s<t$,
\bea
\label{roughpath}
\underline B_{s,t} := \Big[\int_s^t B^i_{s, r}\circ dB^j_r\Big]_{1\le i, j\le d}; \q A_{s,t} := \underline B_{s,t}  -( \underline B_{s,t} )^T, 
\eea
then we can write
\bea
\label{R-rough}
\left\{\ba{lll}
u_{t_\d}= u_t +   \pa_t u_t \d +  \pa_{\o} u_t \cd B_{t,t_\d} +  \pa_{\o\o}^2 u_t : \underline B_{t, t_\d}    + R_2(t,\d);\\
 u_{t^-_\d} = u_t - \pa_t u_{t} \d -   \pa_{\o} u_t \cd B_{t^-_\d,t}   + \pa^2_{\o\o} u_t :  \underline B_{t^-_\d, t}  + R_2(t, -\d).
\ea\right.
\eea  
It is worth noting that $\underline B_{s,t}$ and $A_{s,t}$ are essentially the ``{\it Step-2 signature}" and the ``{\it   L\'evy area}",
respectively,  in Rough Path theory (cf. e.g. \cite{FV}). 
%{\color{red} Indeed, the results of this paper can be easily translated into the rough path language. We leave the details to
% interested readers.} 

(ii) Note that 
\beaa
\int_s^t B^i_{s, r}\circ dB^j_r + \int_s^t B^j_{s, r}\circ dB^i_r = B^i_{s, t}B^j_{s,t},&\mbox{or equivalently},&\underline B_{s,t}  + ( \underline B_{s,t} )^T = B_{s,t} B_{s,t}^T.
\eeaa 
Then (\ref{R-rough}) becomes
\bea
\label{R-rough2}
\left.\ba{lll}
\dis u_{t_\d}= u_t +   \pa_t u_t \d +  \pa_{\o} u_t \cd B_{t,t_\d} + {1\over 2} \pa_{\o\o}^2 u_t :  B_{t, t_\d} B_{t,t_\d}^T +  {1\over 2} \pa_{\o\o}^2 u_t : A_{t,t_\d}  + R_2(t,\d);\ms\\
\dis u_{t^-_\d} = u_t - \pa_t u_{t} \d -   \pa_{\o} u_t \cd B_{t^-_\d,t}   + {1\over 2} \pa_{\o\o}^2 u_t :  B_{t^-_\d,t} B_{t^-_\d,t}^T +  {1\over 2} \pa_{\o\o}^2 u_t : A_{t^-_\d, t}   + R_2(t, -\d).
 \ea\right.
\eea  
Clearly, if  $\pa^2_{\o\o} u$ is symmetric, in particular when $u(t,\o) = v(t,\o_t)$ for some deterministic smooth function $v$, we have $\pa_{\o\o}^2 u_t : A_{t,t_\d}  = \pa_{\o\o}^2 u_t : A_{t^-_\d, t} = 0$, and thus
\bea
\label{R-Markov}
\left.\ba{lll}
\dis u_{t_\d}= u_t +   \pa_t u_t \d +  \pa_{\o} u_t \cd B_{t,t_\d} + {1\over 2} \pa_{\o\o}^2 u_t :  B_{t, t_\d} B_{t,t_\d}^T   + R_2(t,\d);\ms\\
\dis u_{t^-_\d} = u_t - \pa_t u_{t} \d -   \pa_{\o} u_t \cd B_{t^-_\d,t}   + {1\over 2} \pa_{\o\o}^2 u_t :  B_{t^-_\d,t} B_{t^-_\d,t}^T    + R_2(t, -\d).
 \ea\right.
\eea  
This is exactly the standard Taylor expansion. We shall emphasize though, in general  $\pa^2_{\o\o} u$ is not symmetric (see
Example \ref{eg-noncommute}-(ii)), thus the Taylor expansion (\ref{R-rough2})  should have a correction term ${1\over 2} \pa_{\o\o}^2 u : A$.
 \qed}
\end{rem}

 \begin{rem}
 \label{rem-Karandikar}
 {\rm By Bichteler \cite{Bichteler} or Karandikar \cite{Karandikar}, one may interpret $u_t d B^i_t$ in a pathwise manner, whenever $u$ is continuous in $t$. In particular, $\underline B_{s,t}$ and $A_{s,t}$ can be understood pathwisely.
 \qed}
 \end{rem}
As we pointed out in the Introduction, the main results of this paper are the (pathwise) remainder estimates. Since the proof
of the second order estimate is similar to that of the inductional argument for the $m$-th order estimate, we shall prove a general
result directly.
%easier than the general result, we shall prove the direct for It\^o processes. 

\section{Taylor Expansion for It\^{o} Processes (General Case)}
\label{sect-temp-proof}
\setcounter{equation}{0}

We now consider the general form of pathwise Taylor expansion up to any order $m$. 
 Denote, for $0\le t_1<t_2$ and $\e>0$,
\bea
\label{D}
D := \Big\{(t,\d)\in [0, \infty) \times \dbR\backslash \{0\}:  t+\d\ge 0 \Big\}, D^\e_{[t_1, t_2]} := \{(t,\d) \in D: t_1\le t\le t_2, |\d|\le \e\}.
\eea
For any $m\ge 0 $ and $u \in \cH^{[m]}_2(\L)$,  in light of (\ref{Taylor+2}) and (\ref{Taylor-2}) we shall  define the $m$-th order remainder  by: for any $(t,\d)\in D$ and $\o\in \O$,
\bea
\label{Taylorm}
u(t+\d,\o)  =  \sum_{|\th|\le m} \cD^\th_\o u(t,\o) \cI^\th_{t, t+\d} +R_m(u,t,\o,\d).
\eea 
We emphasize that $\d$ can be negative here, and the right side of (\ref{Taylorm})  is pathwise, in light of Remark \ref{rem-Karandikar}. Moreover, when there is no confusion, we shall always omit the variable $\o$.
 
 The main result of this section is the following  pathwise estimate for the remainder $R_m$. 
 %whose proof is deferred to next section.
\begin{thm}
\label{thm-high} 
Assume that $u\in\cH^{[m+2]}_{p_0}(\L)$ for some $m\ge 0$ and $p_0>  2$. Then for any %$2<p< p_0$, and any $0<\a<1$ such that $\frac2{1-\a}<p$, 
$0<\a<1-{2\over p_0}$ and $ p < p_0$, it holds that, for any  $T>0$, 
\bea
\label{highRest}
\dbE\Big\{ \sup_{(t,\d) \in D^1_{[0,T]}}  \Big|{R_m(u, t; \d)\over |\d|^{m+\a\over 2}}\Big|^p\Big\} <\infty.
\eea
\end{thm}

To prove Theorem \ref{thm-high} we need the  following crucial estimate. Since its proof is quite lengthy, we shall 
 complete its proof after we prove Theorem \ref{thm-high}. % \ref{sect-temp-proof} below.
\begin{prop}
\label{prop-high} 
Assume that $u\in\cH^{[m+2]}_{p_0}(\L)$ for some $m\ge 0$ and $p_0>  2$. Then for any $p< p_0$, $t_0\ge 0$, and $\e>0$, it holds that
   \bea
  \label{temp-est}
  \dbE\Big[\sup_{(t, \d) \in D_{[t_0, t_0+\e]}^\e}  |R_m(u, t; \d )|^{p} \Big]  \le C\e^{p(m+1)\over 2},
  \eea
  where $C$ may depend on $\|u\|_{m+2, p_0, T}$ for  some $T\ge t_0+2\e$.
  \qed
\end{prop}

\no [{\it Proof of Theorem \ref{thm-high}.}] 
 In what follows we shall fix $T$, and allow the generic constant $C>0$ to depend on $\|u\|_{m+2, p_0, T+1}$.  Clearly it suffices to prove \reff{highRest} for large $p$, and we thus assume without loss of generality that ${2\over 1-\a} < p < p_0$.

For any $0<\e<1$, set $t_i := i\e$, $i= 0,\cds, [{T\over \e}]+1$. Then, by Proposition \ref{prop-high} we have
 \beaa
%\label{Rme}
 \dbE\Big[\sup_{(t,\d)\in D^\e_{[0, T]}}  |R_m(u,t,\d)|^{p}\Big] \le \sum_{i=0}^{[{T\over \e}]}  \dbE\Big[\sup_{(t,\d)\in D^\e_{[ t_i, t_i+\e]}} |R_m(u,t,\d)|^{p}\Big] \le C\e^{{p(m+1)\over 2} -1}.
\eeaa
Consequently, %for any $p\in (2, p_0)$, $\a\in (0,1)$ such that ${2 \over 1-\a} < p < p_0$, 
 since  $0<\a<1-{2\over p_0}, {2\over 1-\a} < p < p_0$, it holds that
\beaa
&&\dbE\Big[\sup_{(t,\d)\in D^1_{[0,T]}} \Big|{R_m(u,t,\d) \over \d^{m+\a\over 2}}\Big|^p\Big]
%&=&\dbE\Big[\sup_{n\ge 0} \sup_{(t,\d): 0\le t\le T, 2^{-(n+1)}\le |\d| \le 2^{-n}, t+\d\ge 0}  \Big|{R_m(u,t,\d) \over \d^{m+\a\over 2}}\Big|^p\Big]\\
\le \sum_{n=0}^\infty \dbE\Big[\sup_{(t,\d)\in D^1_{[0,T]}:  2^{-(n+1)}\le |\d| \le 2^{-n}}  \Big|{R_m(u,t,\d) \over \d^{m+\a\over 2}}\Big|^p\Big]\\
&\le& \sum_{n=0}^\infty2^{p(m+\a)(n+1)\over 2} \dbE\Big[\sup_{(t,\d)\in D^{2^{-n}}_{[0,T]}}  |R_m(u,t,\d)|^p\Big]\le C \sum_{n=0}^\infty2^{p(m+\a)(n+1)\over 2}  2^{-n({p(m+1)\over 2} -1)} \\
&=& C2^{p(m+\a)\over 2} \sum_{n=0}^\infty 2^{-{n\over 2}[p(1-\a)-2]}<\infty,
\eeaa
completing the proof.
\qed

\ms

\begin{rem}
\label{rem-high}
{\rm (i) The estimate (\ref{highRest}) amounts to saying that, for each $m\ge 0$, $T>0$, there exist a set $\O_m\subseteq \O$ with $\hP(\O_m)=1$ and  a nonnegative random variable $C_{m,T}$, such that 
\beaa
|R_m(u, t, \o; \d )|\le C_{m,T}(\o) |\d|^{m+\a\over 2}, ~~\forall \o\in \O_m, \forall (t,\d) \in D^1_{[0,T]}.
\eeaa
In Theorem \ref{thm-high-space} below, this pathwise estimate also holds true locally uniformly in the spatial variable $x$.

(ii) We should point out that in (\ref{Taylorm}) $\d<0$ is allowed. That is, the temporal expansion can be ``backward". 
Such an expansion, along with the pathwise estimates, is 
%Another main feature of our result is the backward expansion.Both this feature and the above pathwise estimate will be 
crucial for the study of viscosity solutions of SPDEs in \cite{BMZ}. In our previous works  \cite{BM3,BBM} 
%studied both forward and backward expansion 
these results were also obtained
in the case $m=2$, $d=1$, but the present treatment is  more direct and the conditions are easier to verify. (See \S\ref{sect-BBM} 
for a  detailed comparison with the result in \cite{BBM}.) 
%In particular, \cite{BBM} obtained pathwise estimates for $R_m$, which is closest to the present paper.  While our result is more general (higher order expansion in multiple dimensions), our path derivative provides an intuitive meaning of the somewhat mysterious coefficients used in \cite{BM3} and \cite{BBM}, which in particular enables us to transform the SPDE  to a path dependent PDE in a natural way. See more details in Section .

(iii) There have been many works on stochastic Taylor expansions (see, e.g., the books of  \cite{FV} and \cite{KloPla}, and the references
cited therein). We also note that in a recent work \cite{LO}, the Dupire-type path-derivatives were also used. The main difference between
the existing results and ours, however, lies in that in these works the remainder $R_m$ is estimated in $L^2$-sense or in probability, which is not desirable for our study of viscosity solutions. Moreover, no backward expansion was considered in these works.
\qed}
\end{rem}

%\section\no{\bf Proof of Proposition \ref{prop-high}.} ~
In the rest of this section we prove the key estimate (\ref{temp-est}) in Proposition \ref{prop-high}. To simplify the presentation, we split the 
proof into several lemmas that are interesting in their own rights. We begin by establishing a representation formula for $R_m$, extending (\ref{R+2}) and (\ref{R-2}).  In what follows we denote $t_\d := t+\d$ and $t^-_\d := t-\d$, for $\d>0$.
\begin{lem}
\label{lem-f-rep}
Let $u\in \cH^{[m+2]}_{2}(\L)$ for some $m\ge 0$. Then for any $\d>0$, it holds that: 
\bea
\label{f-rep}
R_{m}(u, t, \d) &=& \sum_{|\th| = m+1} \cI^\th_{t, t_\d}( \cD^\th_\o u) +  \sum_{|\th| = m} \cI^\th_{t, t_\d}\Big(\int_t^\cd \pa_t \cD^\th_\o u_s ds\Big).
\eea
Furthermore, denoting $\th = (\th_1, \tilde\th)$, then for $t\ge \d$ one has
\bea
\label{b-rep}
R_{m}(u, t, -\d) =\!\!\! \sum_{|\th| = m+1}  (-1)^{|\th|_0} \int_{t^-_\d}^t \Big(\cD^\th_\o u_s \; \cI^{-\tilde \th}_{t^-_\d, s} \Big) d_{\th_1}s -\!\!\! \sum_{|\th| = m} (-1)^{|\th|_0}\int_{t_\d^-}^t  \Big[\pa_t \cD^\th_\o u_s\; \cI^{-\th}_{t^-_\d, s} \Big]ds.
\eea
\end{lem}

\no{\it Proof.}  (i) We  first verify (\ref{f-rep}) by induction. For $m=0$ we recall the notational convention (\ref{th=0}). Then it is
readily seen that the right side of (\ref{f-rep}) reads
\beaa
 \sum_{|\th| = 1} \int_t^{t_\d} \cD^\th_\o u_s d_\th s +  \int_t^{t_\d} \pa_t  u_sds = \sum_{i=1}^d \int_t^{t_\d} \pa_{\o^i} u_s \circ dB^i_s +  \int_t^{t_\d} \pa_t  u_sds.
 \eeaa
  Thus the equality follows immediately from the functional It\^{o} formula (\ref{Ito2}).

Now assume (\ref{f-rep}) holds for $m$. Then 
\beaa
R_{m+1}(u,t,\d)& =& R_m(u,t,\d) -  \sum_{|\th| = m+1}  \cD^\th_\o u_t \cI^\th_{t, t_\d}\\
&=&\sum_{|\th| = m+1} \cI^\th_{t, t_\d}( \cD^\th_\o u) +  \sum_{|\th| = m} \cI^\th_{t, t_\d}\Big(\int_t^\cd \pa_t \cD^\th_\o u_s ds\Big) -  \sum_{|\th| = m+1} \cD^\th_\o u_t \cI^\th_{t, t_\d}\\
&=&  \sum_{|\th| = m+1}\cI^\th_{t, t_\d}( [\cD^\th_\o u]_{t,\cd})  +   \sum_{|\th| = m} \cI^\th_{t, t_\d}\Big(\int_t^\cd \pa_t \cD^\th_\o u_s ds\Big).
\eeaa
%where the last equality used the fact that $\d>0$ and thus we may move the term $ \cD^\th_\o u (t)$ into the stochastic integral.
Applying the functional It\^{o} formula (\ref{Ito2}) on  $\cD^\th_\o u_s$ we obtain
\beaa
&&R_{m+1}(u,t,\d) \\
&=&  \sum_{|\th| = m+1}\cI^\th_{t, t_\d}\Big( \int_t^\cd \pa_t \cD^\th_\o u_s ds + \sum_{i=1}^d \int_t^\cd \pa_{\o^i} \cD^\th_\o u_s \circ dB^i_s\Big)  +   \sum_{|\th| = m} \cI^\th_{t, t_\d}\Big(\int_t^\cd \pa_t \cD^\th_\o u_s ds\Big)\\
&=&  \sum_{|\th| = m+1}\cI^\th_{t, t_\d}\Big( \int_t^\cd \pa_t \cD^\th_\o u_s ds\Big) \\
&&+\sum_{|\th| = m+1} \sum_{i=1}^d\cI^\th_{t, t_\d}\Big( \int_t^\cd \pa_{\o^i} \cD^\th_\o u_s \circ dB^i_s\Big)  +   \sum_{|\th| = m} \cI^\th_{t, t_\d}\Big(\int_t^\cd \pa_t \cD^\th_\o u_s ds\Big). 
\eeaa
One may check directly that the last line above is exactly equal to $\sum_{|\th| = m+2} \cI^\th_{t, t_\d}(\cD^\th_\o u)$. Namely
(\ref{f-rep}) holds for $m+1$. Thus (\ref{f-rep}) holds for all $m$.

(ii) We now prove (\ref{b-rep}), again by induction. For $m=0$ the argument is similar to (i). Assume now (\ref{b-rep}) holds for $m$. Then
\beaa
&&R_{m+1}(u,t,-\d)  = R_m(u,t,-\d) -  \sum_{|\th| = m+1} (-1)^{|\th|_0} \cD^\th_\o u_t  \;\cI^{-\th}_{t^-_{\d}, t} \\
&&\qq= \sum_{|\th| = m+1}  (-1)^{|\th|_0}\int_{t^-_{\d}}^t\Big[ \cD^\th_\o u_s\;\cI^{-\tilde\th}_{t^-_{\d}, s}\Big] d_{\th_1} s \\
&&\qq\q- \sum_{|\th| = m} (-1)^{|\th|_0}\int_{t^-_{\d}}^t  \Big[\pa_t \cD^\th_\o u_s \;\cI^{-\th}_{t^-_{\d}, s}\Big]ds-  \sum_{|\th| = m+1} (-1)^{|\th|_0} \cD^\th_\o u_t \int_{t^-_{\d}}^t \;\cI^{-\tilde\th}_{t^-_{\d}, s} d_{\th_1}s.
\eeaa
Applying integration by parts formula we have % (under Stratonovich integration), %or equivalently applying the It\^{o} formula 
\beaa
%\label{-th}
&& \cD^\th_\o u_t \int_{t^-_{\d}}^t \;\cI^{-\tilde\th}_{t^-_{\d}, s} d_{\th_1}s -\int_{t^-_{\d}}^t\Big[ \cD^\th_\o u_s\;\cI^{-\tilde\th}_{t^-_{\d}, s}\Big] d_{\th_1} s\\
&=&\int_{t^-_{\d}}^t \;\cI^{-\th}_{t^-_{\d}, s} \circ d(\cD^\th_\o u_s )  = \int_{t^-_{\d}}^t\;\cI^{-\th}_{t^-_{\d}, s}  \big[\pa_t \cD^\th_\o u_s ds + \sum_{i=1}^d \pa_{\o^i} \cD^\th_\o u_s  \circ dB^i_s\big].\nonumber
\eeaa
Consequently we obtain:
\beaa
&&R_{m+1}(u,t,-\d)
%&=& -  \sum_{|\th| = m+1}  (-1)^{|\th|_0}\int_{t^-_{\d}}^t \int_{t^-_{\d}}^s d_{-\th} r  \Big[\pa_t \cD^\th_\o u_s ds + \sum_{i=1}^d \pa_{\o^i} \cD^\th_\o u_s  \circ dB^i_s\Big]- \sum_{|\th| = m} (-1)^{|\th|_0}\int_{t^-_{\d}}^t  \Big[\pa_t \cD^\th_\o u_s \int_{t^-_{\d}}^s d_{-\th} r\Big]ds\\
=  -  \sum_{|\th| = m+1}  (-1)^{|\th|_0}\int_{t^-_{\d}}^t \Big[\pa_t \cD^\th_\o u_s\;\cI^{-\th}_{t^-_{\d}, s}\Big]ds\\
&&\q\qq  -\sum_{|\th| = m+1}  (-1)^{|\th|_0}\int_{t^-_{\d}}^t \;\cI^{-\th}_{t^-_{\d}, s}    \pa_{\o} \cD^\th_\o u_s  \circ dB_s- \sum_{|\th| = m} (-1)^{|\th|_0}\int_{t^-_{\d}}^t  \Big[\pa_t \cD^\th_\o u_s \;\cI^{-\th}_{t^-_{\d}, s}\Big]ds.
\eeaa
One may now check directly that the last line above is exactly equal to, denoting $\th = (\th_1, \tilde\th)$, 
\beaa
 \sum_{|\th| = m+2}  (-1)^{|\th|_0}\int_{t^-_\d}^t\Big[ \cD^\th_\o u_s \;\cI^{-\tilde\th}_{t^-_{\d}, s}\Big] d_{\th_1} s.
 \eeaa
  Thus (\ref{b-rep}) holds for $m+1$, proving (ii), whence the Lemma.
\qed
  
To simplify notation, in what follows we denote, for any  semi-martingale $\f$,  and any $p\ge 1$, $0\le t_1 < t_2$,
\bea
\label{Iphi}
I(\f, p, t_1, t_2) &:=& \Big(\dbE\Big[ \big(\int_{t_1}^{t_2} |\f_s|^2 ds\big)^{p\over 2} + \big(\int_{t_1}^{t_2} |\pa_\o\f_s| ds\big)^p \Big]\Big)^{1\over p}.
\eea
It is clear that $I$ is increasing in $p$, and
\bea
\label{I1}
 I^p(\f,p, t_1, t_2) \le \sup_{t_1\le t\le t_2} \dbE[|\f_t|^p] (t_2-t_1)^{p\over 2} + \sup_{t_1\le t\le t_2} \dbE[|\pa_\o\f_t|^p] (t_2-t_1)^{p}.
\eea
% We now give an important estimate which will be very useful in our future discussion.

In light of the above representations, the following estimate is crucial.
\begin{lem}
\label{lem-Iest}
Let $t_0\ge 0$, $0<\e<1$, $q>p \ge 1$,  and $\f$ be a  semimartingale. Then, for any $|\th|\ge 1$,   there exists  constant $C = C_{p,q, |\th|_0}$ such that
\bea
\label{H-claim}
\dbE\Big\{\sup_{0<\d  \le \e}\sup_{t_0\le t\le t_0+\e} |\cI^\th_{t, t_\d}(\f)|^p\Big\} &\le&
 C \e^{{p (|\th|-1)\over 2}}I^p(\f, q, t_0, t_0+2\e)
\eea
%\label{H-claim2}
%&&\dbE\Big\{\sup_{0<\d  \le \e}\sup_{t_0\le t\le t_0+\e}\big| \int_t^{t_\d} \f_s \cI^{\tilde\th}_{t,s} d_{\th_1}s\big|^p\Big\} \le C \e^{{p (|\th|-1)\over 2}}\Big(I(\f, q, t_0, t_0+2\e)\Big)^{p}.
% C_{q,\tilde q} \e^{{q (|\th|-1)\over 2}}\Big\{ \dbE\Big[ \big(\int_{t_0}^{t_0+2\e} [|\f_s|^2+|\pa_\o \f_s|^2] ds\big)^{\tilde q\over 2} \Big]\Big\}^{q\over \tilde q}.
\end{lem}

\no{\it Proof.}   Let $a(\f, \th, p)$  denote the left side of (\ref{H-claim}), and $I(\f, p) := I(\f, p, t_0, t_0+2\e)$. 
Without loss of generality, we may assume $I(\f, q) <\infty$.
%
%
%First, we note that if the process $\f$ doest not have the desired integrability, then the right hand side of (\ref{H-claim}) is 
%$\infty$, and the inequality becomes trivial. So from this point on we always assume that $\f$ is such that the right hand side of 
%(\ref{H-claim}) is finite.
We proceed  by induction on $n:=|\th|_0$. %Let {\color{red}$a(\f, \th, p)$} denote the left side of (\ref{H-claim}) and $t_\d := t+\d$. 

(i) First assume $n=1$, namely $\th = (\th_1)$.   We estimate $a(\f, \th, p)$ in two cases.

{\it Case 1.} $|\th|=2$, namely $\th_1=0$. Then $\cI^\th_{t, t_\d}(\f) = \int_t^{t_\d} \f_s ds$, and thus
  \beaa
a(\f,\th, p) \le \dbE\Big[ \big(\int_{t_0}^{t_0+2\e}  |\f_s| ds\big)^p\Big]\le C\e^{p\over 2} I(\f, p)^p\le  C\e^{p\over 2} I(\f, q)^p.
%\dbE\Big[ \big(\int_{t_0}^{t_0+2\e}  |\f_s|^2 ds\big)^{p\over 2}\Big],
 \eeaa
% and the result is  obvious. 
 
 \ms
 
 {\it Case 2.} $|\th|=1$, namely $\th_1 = i$ for some  $i=1,\cds, d$. Then
 \beaa
\cI^\th_{t, t_\d}(\f) &=& [\int_{t_0}^{t_\d} - \int_{t_0}^t]  \f(s)  dB^i_s +{1\over 2} \int_t^{t_\d}  \pa_{\o^i}\f_s ds,
\eeaa
and thus
\beaa
a(\f,\th,p) 
 &\le& C \dbE\Big[ \sup_{t_0\le t\le t_0+2\e} \big|\int_{t_0}^{t}   \f_s  dB^i_s \big|^p \Big] + C\dbE\Big[ \big(\int_{t_0}^{t_0+2\e}  |\pa_{\o^i}\f_s| ds\big)^p\Big],
 \eeaa
which, together with  the  Burkholder-Davis-Gundy inequality, implies  \reff{H-claim}  immediately.
 
 \ms
 
(ii) We next prove \reff{H-claim} by induction. Assume it  holds true for $ n$ and we now assume  $|\th|_0=n+1$.  Denote $\th = (\th_1, \tilde\th)$,  $\psi_s := \int_{t_0}^s \f_r d_{\th_1} r$, and $\tilde q := {p+q\over 2}$. Notice that $|\tilde \th|_0=n$ and $I(1,p) = \sqrt{2\e}$ for any $p$, then we may use the induction assumption  and obtain
\beaa
a(\f, \th, p) &=& \dbE\Big[\sup_{0<\d  \le \e}\sup_{t_0\le t\le t_0+\e} \big|\cI^{\tilde\th}_{t,t_\d}(\psi) - \psi_t \cI^{\tilde\th}_{t,t_\d}(1) \big|^p\Big] \\
&\le& Ca(\psi, \tilde \th, p) + C\Big(a(\f, (\th_1), \tilde q)\Big)^{p\over \tilde q}\Big( a(1, \tilde\th, {2p \tilde q\over q-p})\Big)^{q-p\over 2\tilde q}\\
%&\le& \dbE\Big[\sup_{0<\d  \le \e}\sup_{t_0\le t\le t_0+\e} \big|\cI^{\tilde\th}_{t,t_\d}(\psi) \big|^p\Big] \\
%&&+ \Big(\dbE\Big[\sup_{t_0\le t\le t_0+\e} |\psi_t|^{p+q\over 2}\Big]\Big)^{2p\over p+q} \Big( \dbE\Big[\sup_{0<\d  \le \e}\sup_{t_0\le t\le t_0+\e} | \cI^{\tilde\th}_{t,t_\d}(1) |^{p(p+q)\over q-p}\Big]\Big)^{q-p\over p+q}\\
&\le& C\e^{p(|\tilde\th|-1)\over 2}  \big( I(\psi, \tilde q)\big)^p  +C \e^{p(|(\th_1)|-1)\over 2}\big( I(\f, q)\big)^{p} \e^{p(|\tilde \th|-1)\over 2} \big(I(1,{3p \tilde q \over q-p})\big)^{p} \\
&=&C\e^{p(|\tilde\th|-1)\over 2}  \big( I(\psi, \tilde q)\big)^p  +C \e^{p(|\th|-1)\over 2}\big( I(\f, q)\big)^{p}.
 \eeaa
% where the last inequality used the obvious fact that $I(1,p) = \sqrt{2\e}$ for any $p$. 
Then clearly \reff{H-claim} with $|\th|_0=n+1$ follows from the following claim:
  \bea
 \label{H-claim3}
  I(\psi, \tilde q) &\le& C \e^{|(\th_1)|\over 2} I(\f,  q).
  \eea
 We again proceed in two cases.

{\it Case 1.} $|(\th_1)|=2$, namely $\th_1=0$.  Then  $\psi_s = \int_{t_0}^s \f_r  dr$ and $\pa_{\o^j} \psi_s = 0$. Thus, 
\beaa
I(\psi, \tilde q)^{\tilde q}&\le&\dbE\Big[ \big(\int_{t_0}^{t_0+2\e} |\int_{t_0}^s \f_r dr|^2 ds\big)^{\tilde q\over 2}\Big] \le \e^{\tilde q\over 2} \dbE\Big[ \big(\int_{t_0}^{t_0+2\e}  |\f_r| dr\big)^{\tilde q}\Big]  \le  \e^{\tilde q} I(\f,  q)^{\tilde q}.
\eeaa

{\it Case 2.} $|(\th_1)|=1$, namely $\th_1 = i$ for some $i=1,\cds, d$. Then
\beaa
\psi_s = \int_{t_0}^s \f_r \circ dB^i_r =  \int_{t_0}^s \f_r dB^i_r  + {1\over 2} \int_{t_0}^s \pa_{\o^i} \f_r dr,\q \pa_{\o^j} \psi_s = \f_s 1_{\{j=i\}},
\eeaa
and thus 
\beaa
I(\psi, \tilde q)^{\tilde q} &\le&C \dbE\Big[ \big(\int_{t_0}^{t_0+2\e} [|\int_{t_0}^s \f_r dB^i_r|^2+|\int_{t_0}^s \pa_{\o^i} \f_r dr|^2] ds\big)^{\tilde q\over 2} + \big(\int_{t_0}^{t_0+2\e} |\f_s|ds\big)^{\tilde q} \Big]\\
&\le& C\e^{\tilde q\over 2} \dbE\Big[ \sup_{t_0\le  s \le t_0+2\e} \big[|\int_{t_0}^s \f_r dB^i_r  |^{\tilde q}\big]+\big[\int_{t_0}^{t_0+2\e}  |\pa_{\o^i} \f_r| dr\big]^{\tilde q} + \big(\int_{t_0}^{t_0+2\e} |\f_s|^2ds\big)^{\tilde q\over 2} \Big].
\eeaa
Then \reff{H-claim3} follows immediately from  the Burkholder-Davis-Gundy inequality.
 \qed

We are now ready to prove Proposition \ref{prop-high}. 

\no[{\it Proof of Proposition \ref{prop-high}}.]   First,  by \reff{f-rep} we have, for $\d>0$,
\beaa
R_{m}(u, t, \d) &=& \sum_{|\th| = m+1} \cI^\th_{t, t_\d}( \cD^\th_\o u) +  \sum_{|\th| = m}\Big[ \cI^\th_{t, t_\d}\Big(\int_{t_0}^\cd \pa_t \cD^\th_\o u_s ds\Big) -  \cI^\th_{t, t_\d}(1)\int_{t_0}^t \pa_t \cD^\th_\o u_s ds \Big].
\eeaa
Denote $D^\e_+ := \{(t,\d): 0<\d\le \e, t_0\le t\le t_0+\e\}$ and note that $\pa_\o\int_{t_0}^t \pa_t \cD^\th_\o u_s ds = 0$. Then, combine Lemma \ref{lem-Iest} and \reff{I1} and recall $T\ge t_0+2\e$, we have
\bea
\label{Rmest1}
&& \dbE\Big[\sup_{(t,\d)\in D^\e_+} |R_m(u, t, \d )|^{p}\Big]\\
 &\le& C \sum_{|\th|=m+1} \dbE\Big[\sup_{(t,\d)\in D^\e_+} |\cI^\th_{t, t_\d}( \cD^\th_\o u)|^p\Big] +C \sum_{|\th|=m}   \dbE\Big[\sup_{(t,\d)\in D^\e_+} |\cI^\th_{t, t_\d}( \int_{t_0}^\cd \pa_t \cD^\th_\o u_s ds)|^p\Big]\nonumber\\
 &&+ C \sum_{|\th| = m}\Big(\dbE\Big[\big(\int_{t_0}^{t_0+2\e} |\pa_t \cD^\th_\o u_s| ds\big)^{p_0}\Big]\Big)^{p\over p_0} \Big(\dbE\Big[\big(\sup_{(t,\d)\in D^\e_+}|\cI^\th_{t, t_\d}(1)|\big)^{pp_0\over p_0-p}\Big]\Big)^{p_0-p\over p_0} \nonumber\\
 &\le& C \sum_{|\th|=m+1}\Big[ \| \cD^\th_\o u\|_{p,T}^p \e^{p(m+1)\over 2} + \|\pa_\o \cD^\th_\o u\|_{p,T}^p \e^{p(m+2)\over 2}\Big]\nonumber\\
 &&+C\e^{pm\over 2} \sum_{|\th|=m}  \sup_{t_0\le t\le t_0+2\e} \dbE\Big[\Big|\int_{t_0}^t \pa_t \cD^\th_\o u_s ds\Big|^p\Big] + C \sum_{|\th| = m}\Big(\|\pa_t \cD^\th_\o u\|_{p_0, T}^{p} \e^{p}\Big) \e^{pm\over 2} \nonumber\\
&\le& C \sum_{|\th|=m+1}\Big[ \| \cD^\th_\o u\|_{p_0,T}^p \e^{p(m+1)\over 2} + \|\pa_\o \cD^\th_\o u\|_{p_0,T}^p \e^{p(m+2)\over 2}\Big]\nonumber\\
&&+C \sum_{|\th|=m} \Big[ \|\pa_t \cD^\th_\o u\|_{p_0,T}^p \e^{p(m+2)\over 2} + \|\pa_t \cD^\th_\o u\|_{p_0, T}^{p} \e^{p(m+1)\over 2} \Big] \nonumber\\
&\le& C\|u\|_{m+2, p_0, T} \e^{p(m+1)\over 2}.\nonumber
 \eea

Next, for $\th = (\th_1, \tilde \th) = (\th_1, \th_2,\cds, \th_n)$,  applying the integration by parts formula and recalling \reff{DI2} and \reff{-th}, 
we have, for $\f := \cD^\th_\o u$,
\beaa
\int_{t_\d^-}^t \big(\f_s \cI^{-\tilde\th}_{t_\d^-,s}\big) d_{\th_1} s &=& \Big(\int_{t_\d^-}^t \f_s  d_{\th_1} s\Big) \cI^{-\tilde\th}_{t_\d^-,t} - \int_{t_\d^-}^t \big(\int_{t_\d^-}^s \f_r  d_{\th_1} r\big) \circ d\cI^{-\tilde\th}_{t_\d^-,s}\\
&=&\cI^{(\th_1)}_{t^-_\d, t}(\f) \cI^{(\th_n,\cds, \th_2)}_{t^-_\d, t}(1) - \int_{t_\d^-}^t \Big(\cI^{(\th_1)}_{t^-_\d, s}(\f) \cI^{(\th_n,\cds, \th_3)}_{t^-_\d, s}\Big) d_{\th_2}s
\eeaa
Repeating the above arguments we obtain
\bea
\label{Rm-rep}
\int_{t_\d^-}^t \big(\f_s \cI^{-\tilde\th}_{t_\d^-,s}\big) d_{\th_1} s = \sum_{i=1}^n (-1)^{i-1} \cI^{(\th_1,\cds, \th_i)}_{t^-_\d, t}(\f) \cI^{(\th_n,\cds, \th_{i+1})}_{t^-_\d, t}(1).
\eea
Then, by changing variable $t-\d$ to $t$ and denoting $q:= {p+p_0\over 2}$, we have
\bea
\label{Rmest2}
&& \dbE\Big[\sup_{0<\d  \le \e}\sup_{t_0\vee \d \le t\le t_0+\e} |\int_{t_\d^-}^t \big(\cD^\th_\o u_s \cI^{-\tilde\th}_{t_\d^-,s}\big) d_{\th_1} s|^{p}\Big]\\
&\le& C\sum_{i=1}^n \dbE\Big[\sup_{0<\d  \le \e}\sup_{(t_0-\e)^+\le t\le t_0+\e} | \cI^{(\th_1,\cds, \th_i)}_{t, t_\d}(\cD^\th_\o u) \cI^{(\th_n,\cds, \th_{i+1})}_{t,t_\d}(1)
|^{p}\Big]\nonumber\\
&\le&  C\sum_{i=1}^n \Big( \dbE\Big[\sup_{0<\d  \le \e}\sup_{(t_0-\e)^+\le t\le t_0+\e} | \cI^{(\th_1,\cds, \th_i)}_{t, t_\d}(\cD^\th_\o u)|^{q}\Big]\Big)^{p\over q}\times\nonumber\\
&& \Big( \dbE\Big[\sup_{0<\d  \le \e}\sup_{(t_0-\e)^+\le t\le t_0+\e} |\cI^{(\th_n,\cds, \th_{i+1})}_{t,t_\d}(1)|^{2pq \over p_0-p}\Big]\Big)^{p_0-p\over 2q}\nonumber\\
&\le& C\sum_{i=1}^n  \e^{p(|(\th_1,\cds,\th_i)|-1)\over 2} I^p(\cD^\th_\o u, , p_0, (t_0-\e)^+, t_0+2\e) \e^{p|(\th_n,\cds, \th_{i+1})|\over 2}\nonumber\\
&=&C \e^{p(|\th|-1)\over 2} I^p(\cD^\th_\o u, , p_0, (t_0-\e)^+, t_0+2\e),\nonumber
 \eea
 thanks to Lemma  \ref{lem-Iest}. Now following similar arguments as in \reff{Rmest1}, one may easily derive from \reff{b-rep} that
 \bea
\label{Rmest3}
\dbE\Big[\sup_{0<\d  \le \e}\sup_{t_0\vee \d \le t\le t_0+\e} |R_m(u, t, -\d )|^{p}\Big]&\le& C\|u\|_{m+2, p_0, T} \e^{p(m+1)\over 2}.
 \eea
 Finally, note that
 \beaa
\sup_{(t,\d) \in D^\e_{[t_0, t_0+\e]}}  |R_m(u, t; \d )|&\le& \sup_{0<\d  \le \e}\sup_{t_0\le t\le t_0+\e} |R_m(u, t; \d )|+\sup_{0<\d  \le \e}\sup_{t_0\vee \d \le t\le t_0+\e} |R_m(u, t; -\d )|.
   \eeaa
Combining \reff{Rmest1} and \reff{Rmest3}  we obtain \reff{temp-est}.
\qed

\section{Pathwise Taylor Expansion for Random Fields}
\label{sect-space}
\setcounter{equation}{0}

In this  section we extend our results to It\^o random fields.  %and then we weaken the regularity assumption on $u$.  
Again, denote $t_\d := t+\d$ and we emphasize that $\d$ could be negative.
%\subsection{Pathwise Taylor expansion}
Let us denote
\bea
\label{D-space}
\left.\ba{lll}
\dis \hat D := \Big\{(t,x,\d, h)\in [0, \infty) \times \cO \times \dbR\times \dbR^{d'}:  (t_\d, x+h) \in Q\Big\};\\
\dis \hat D_N := \Big\{ (t,x,\d, h):  (t,x)\in Q_N,  |\d|+|h|^2 \le {1\over N(N+1)}, t_\d\ge 0 \Big\}.
\ea\right.
\eea
We remark that $(t_\d, x+h)\in Q_{N+1}$ for any $(t,x,\d, h) \in \hat D_N$.
Furthermore, for any $m\ge 0 $ and $u \in \cH^{[m]}_2(\hat\L)$,  in light of (\ref{Taylorm}), and noting that the spatial derivative
$\pa_x$ commutes with all the path derivatives, we shall  define the $m$-th order Taylor expansion  by: 
\bea
\label{Taylorm-space}
u(t+\d, x+h)  =  \sum_{|(\th,\ell)|\le m}{1\over \ell! } \big(\cD^\ell_x \cD^\th_\o u\big)(t,x)  h^{\ell}\cI^\th_{t, t_\d}  +R_m(u,t,x,\d, h),
\eea
for any $(t,x,\d,h)\in D$ and $\o\in \O$.  %We emphasize again that  $\d$ could be negative.
Clearly,  if $\d=0$, then we recover the standard Taylor expansion in $x$; and if $h={\bf 0}$, then we recover the Taylor expansion (\ref{Taylorm}) for It\^o processes. Moreover, if $m=2$, then we have
\bea
\label{Taylor2-space}
u(t+\d, x+h)  &=&  \Big[u + \pa_t u \d + \pa_x u \cd h + \pa_\o u \cd B_{t, t_\d} + {1\over 2} \pa^2_{xx} u : hh^T \\
&& + \pa_{x\o} u : h[B_{t,t_\d}]^T + \pa^2_{\o\o} u : \underline B_{t, t_\d}\Big](t,x)   +R_2(u,t,x, \d, h).\nonumber
\eea

Again, we begin with the following simple recursive relations for the remainders.
 \begin{lem}
\label{lem-x-rep}
Let $u \in \cH^{[m+1]}(\hat\L)$ for some  $m\ge 1$. Then, 
\bea
\label{x-rep}
R_{m}(u, t, x,\d, h) =R_m(u,t,x,\d, {\bf 0}) +   \sum_{i=1}^{d'}h_i  \int_0^1 R_{m-1}(\pa_{x_i} u, t, x,\d,  h^i(\k))d\k,
\eea
where   $h^i(\k) := (h_1,\cds, h_{i-1}, \k h_i, 0,\cds, 0)$.
\end{lem}
\no {\it Proof.} 
Given $(t,x,\d, h)\in D$, we write
 \bea
\label{R}
u(t_\d,x+ h)  -  u(t,  x)= E_0+ \sum_{i=1}^{d'} E_i,
\eea
where $E_0 := u(t_\d,  x)-  u(t,  x)$ and $E_i :=u(t_\d,x+ h^i(1) ) - u(t_\d,x + h^i(0))$. Note that by applying the temporal expansion 
we have
\bea
\label{I2}
E_0&=&\sum_{1\le |\th| \le m}\cD^\th_\o u (t,x)\;\cI^\th_{t, t_\d} + R_m(u,t,x,\d, {\bf 0}).
\nonumber\\
&=&\sum_{1\le |(\th, \ell)| \le m, |\ell|=0}{1\over \ell! } \cD^{\ell}_x \cD^\th_\o u (t,x) h^\ell \;\cI^\th_{t, t_\d} + R_m(u,t,x,\d, {\bf 0}).
\eea
On the other hand,  for $i=1,\cds, d'$, using the Taylor expansion for $\pa_xu$ we can write 
\beaa
E_i \!\!&=&\!\! h_i  \int_0^1 \pa_{x_i} u(t_\d, x+ h^i(\k)) d\k \\  
\!\!&= &\!\!h_i   \int_0^1\Big[ \sum_{|(\th,\ell)|\le m-1}{1\over \ell!} \cD^\ell_x\cD^\th_\o (\pa_{x_i}u)(t,x) ( h^i(\k))^{\ell} \;\cI^\th_{t, t_\d}+ R_{m-1}(\pa_{x_i} u, t,x,\d,  h^i(\k))\Big] d\k\\
%.
%\eeaa
%By (\ref{cHn+1}), we have
%\beaa
%E_i
&=&   \sum_{|(\th,\ell)|\le m-1,\ell_{i+1}=\cds=\ell_{d'}=0}{ h_1^{\ell_1}\cds h_{i-1}^{\ell_{i-1}} h_i^{\ell_i+1}\over \ell_1\cds\ell_{i-1}(\ell_i+1)} \pa_{x_i}\cD^\ell_x\cD^\th_\o u(t,x)  \;\cI^\th_{t, t_\d} \\
&&+ h_i \int_0^1 R_{m-1}(\pa_{x_i} u, t,x,\d,  h^i(\k)) d\k\nonumber\\
&=&   \sum_{|(\th,\ell)|\le m, \ell_i\ge 1, \ell_{i+1}=\cds=\ell_{d'}=0}{ h^\ell \over \ell!} \cD^\ell_x\cD^\th_\o u(t,x)  \;\cI^\th_{t, t_\d} + h_i \int_0^1 R_{m-1}(\pa_{x_i} u, t,x,\d,  h^i(\k)) d\k,
\eeaa
where we replaced $\ell = (\ell_1,\cds, \ell_i, 0,\cds, 0)$ with $(\ell_1, \cds, \ell_{i-1}, \ell_i + 1, 0,\cds, 0)$. Then
\beaa
\sum_{i=1}^{d'} E_i =  \sum_{|(\th,\ell)|\le m, |\ell|\ge 1}{ h^\ell \over \ell!} \cD^\ell_x\cD^\th_\o u(t,x) \;\cI^\th_{t, t_\d} +\sum_{i=1}^{d'} h_i \int_0^1 R_{m-1}(\pa_{x_i} u, t,x,\d,  h^i(\k)) d\k.
\eeaa
This, together with (\ref{R}) and (\ref{I2}), implies that
\beaa
u(t_\d,x+ h)  -  u(t,  x) &=& \sum_{1\le |(\th,\ell)|\le m}{ h^\ell \over \ell!} \cD^\ell_x\cD^\th_\o u(t,x)  \;\cI^\th_{t, t_\d}\\
&& + R_m(u,t,x,\d, {\bf 0}) +   \sum_{i=1}^{d'}h_i  \int_0^1 R_{m-1}(\pa_{x_i} u, t, x,\d, h^i(\k))d\k.
\eeaa
Now (\ref{x-rep}) follows immediately from (\ref{Taylorm-space}).
\qed

%\subsection{Estimate of the remainder}
 Our main result  of this section is the following  pathwise estimate for the remainders, extending Theorem \ref{thm-high}.
\begin{thm}
\label{thm-high-space} 
Assume $u\in\cH^{[m+2]}_{p_0}$ for some $m\ge 0$ and $p_0> p_*:= (m+1)d' + 2$. Then for any $p< p_0$ and $0 < \a < 1-{p_*\over p_0}$,
it holds that, for any  $N>0$, 
\bea
\label{highRest-space}
\dbE\Big\{ \sup_{(t,x,\d, h) \in \hat D_N}  \Big|{R_m(u, t,x; \d,  h )\over (|\d| + | h |^2)^{m+\a\over 2}}\Big|^p\Big\} <\infty.
\eea
\end{thm}

\no{\it Proof.}   We  fix $N$ and let $0<\e<{1\over 4 N^2(N+1)^2}$. Then for any $(t,x)\in Q_N$ (recall (\ref{KN})), and $|\d|\le \e, |h|\le \sqrt{\e}$, we have $((t_\d)^+, x+h) \in Q_{N+1}$ and $(t+2\e, x+h)\in Q_{N+1}$. In what follows our generic constant $C$ will depend on $\|u\|_{m+2, p_0, N+1}$. % := \sum_{|(\th, \ell)|\le m+2} \|\cD^\ell_x\cD^\th_\o u\|_{p_0, N+1}$.  
Denote, for $(t_0, x_0)\in Q_N$ and $n\ge 0$,
\bea
\label{DNe}
\hat D_n^{\e}(t_0,x_0):=\{(t,x,\d,h): |\d| + |h|^2\le \e,  t_0\le t \le t_0+\e,   t_\d \ge 0,   |x-x_0|\le \e^{n+1\over 2} \}.
\eea

We split the proof into the following steps. 

(i) We first show that, for any $(t_0,x_0)\in Q_N$ and $p<p_0$,
\bea
\label{R0}
&\hE\Big[\sup_{(t,x,\d, h)\in \hat D_0^{\e}(t_0,x_0)}  |R_0(u,t,x,\d,h)|^p\Big] \le C\e^{p\over 2}.
\eea
Indeed,  note that
$$| R_0(u,t,x,\d,h)|=|u(t_\d, x+ h)-u(t,x)| \le  R_{0,1}+R_{0,2}, 
$$
where $R_{0,1}:=|u(t_\d, x+ h)-u(t_\d, x_0)| +|u(t, x_0)-u(t,x)|$, and $R_{0,2}:=|u(t_\d, x_0)-u(t,x_0)|$.

Note that $R_{0,2}$ is for fixed $x_0$. Applying Proposition \ref{prop-high}  we get
\bea
 \label{R02}
  \dbE\Big[\sup_{(t,x,\d, h)\in \hat D_0^{\e}(t_0,x_0)} |R_{0,2}|^{p}\Big]\le C\e^{p\over 2}.
  \eea
Moreover, note that
  \beaa
   R_{0,1}&=&\Big|(x-x_0+h)\cd \int_0^1 \pa_x u(t _\d, x_0 + \k (x-x_0+h))] d\k\Big| \\
   &&+ \Big| (x-x_0)\cd \int_0^1 \pa_x u(t, x_0 + \k (x-x_0)) d\k\Big|,
  \eeaa
and $x_0 + \k (x-x_0), x_0 + \k (x-x_0+h) \in K_{N+1}$, thanks to \reff{KN2}. % the convexity of $K_{N+1}$ (recall (\ref{KN})). 
Then, since $\pa_x u \in \cH^{[m+1]}_{p_0} \subset \cH^{[0]}_{p_0}$ and $p< p_0$, we have
\beaa
 \dbE\Big[\sup_{(t,x,\d, h)\in D_0^{\e}(t_0,x_0)} |R_{0,1}|^{p}\Big]\le C\e^{p\over 2}.
  \eeaa
This, together with (\ref{R02}), proves (\ref{R0}).  

\ms

(ii) We next show by induction on $m$ that, for any $m\ge 0$, any $(t_0,x_0)\in D_N$ and $p<p_0$,
\bea
\label{Rm}
\hE\Big[\sup_{(t,x,\d, h)\in \hat D_m^{\e}(t_0,x_0)}  |R_m(u,t,x,\d,h)|^p\Big] \le C\e^{p(m+1)\over 2}.
\eea
Indeed, by (\ref{R0}) we have (\ref{Rm}) for $m=0$. Assume (\ref{Rm}) holds true for $m-1$. Applying Lemma \ref{lem-x-rep}, we have
\beaa
&&\hE\Big[\sup_{(t,x,\d, h)\in  \hat D_m^{\e}(t_0,x_0)}  |R_m(u,t,x,\d,h)|^p\Big]  \le C\hE\Big[\sup_{(t,x,\d, h)\in  \hat D_m^{\e}(t_0,x_0)}  |R_m(u,t,x,\d,{\bf 0})|^p\Big] \\
&&\qq  + C\e^{p\over 2} \sum_{i=1}^{d'}\int_0^1 \hE\Big[\sup_{(t,x,\d, h)\in  \hat D_m^{\e}(t_0,x_0)} R_{m-1}(\pa_{x_i} u, t, x,\d, \o, h^i(\k))|^p\Big]d\k.
\eeaa
Note that $\pa_{x_i} u \in \cH^{[m+1]}_{p_0}$ and $(t,x,\d, h^i(\k)) \in   \hat D_m^{\e}(t_0,x_0)\subset  \hat D_{m-1}^{\e}(t_0,x_0)$. Then by induction assumption  we have
\beaa
\hE\Big[\sup_{(t,x,\d, h)\in  \hat D_m^{\e}(t_0,x_0)}  R_{m-1}(\pa_{x_i} u, t, x,\d, \o, h^i(\k))|^p\Big] \le C\e^{pm\over 2}.
\eeaa
So it suffices to prove
\bea
\label{Rm0}
\hE\Big[\sup_{(t,x,\d, h)\in  \hat D_m^{\e}(t_0,x_0)} |R_m(u,t,x,\d,{\bf 0})|^p\Big] \le C\e^{p(m+1)\over 2}.
\eea

To this end, we note that
\bea
\label{Rm1}
R_m(u,t,x,\d,{\bf 0}) = R_m(u,t,x,\d,{\bf 0}) - R_m(u,t,x_0,\d,{\bf 0}) + R_m(u,t,x_0,\d,{\bf 0}). 
\eea
Applying Proposition \ref{prop-high}   again we have
\bea
 \label{Rm2}
  \dbE\Big[\sup_{(t,x,\d, h)\in  \hat D_m^{\e}(t_0,x_0)} |R_m(u,t,x_0,\d,{\bf 0})|^{p}\Big]\le C\e^{p(m+1)\over 2}.
  \eea
On the other hand, recall that
$
R_m(u,t,x,\d,{\bf 0}) = u(t_\d, x) - \sum_{|\th|\le m} \cD^\th_\o u(t,x) \;\cI^\th_{t, t_\d}.
$
Then
\beaa
&& R_m(u,t,x,\d,{\bf 0}) - R_m(u,t,x_0,\d,{\bf 0})\\
 &=& (x-x_0) \cd \int_0^1\Big[ \pa_x u (t_\d, x_0+\k (x-x_0)) - \sum_{|\th|\le m} \pa_x \cD^\th_\o u(t,x_0+\k (x-x_0)) \;\cI^\th_{t, t_\d}\Big]d\k.
\eeaa
Notice that $\pa_x u, \pa_x \cD^\th_\o \in \cH^{[1]}_{p_0}$ for $|\th|\le m$, and $|x-x_0|\le \e^{m+1\over 2}$ for  $(t,x,\d, h)\in  \hat D_m^{\e}(t_0,x_0)$. Then, by Lemma \ref{lem-Iest}, 
\beaa
 \dbE\Big[\sup_{(t,x,\d, h)\in  \hat D_m^{\e}(t_0,x_0)} |R_m(u,t,x,\d,{\bf 0}) - R_m(u,t,x_0,\d,{\bf 0})|^{p}\Big]\le  C\e^{p(m+1)\over 2}.
\eeaa
Plugging this and (\ref{Rm}) into (\ref{Rm1}) we obtain (\ref{Rm0}), which in turn implies (\ref{Rm}).

\ms

(iii) We now claim that
\bea
\label{Rme}
\hE\Big[\sup_{(t,x,\d, h)\in \hat D_N, |\d|+|h|^2\le \e}  |R_m(u,t,x,\d,h)|^p\Big] \le C \e^{{(p-d')(m+1)\over 2} -1}.
\eea
Indeed, set $t_i := i\e$, $i= 0,\cds, [{N\over \e}]+1$, and let $x_j\in K_N$, $j=1,\cds, [(2N\e^{-{m+1\over 2}})^{d'}]+1$ be discrete grids such that the union of their $\e^{m+1\over 2}$-neighborhood covers $K_N$.   Then we have
\beaa
\Big\{(t,x,\d, h)\in \hat D_N:  |\d|+|h|^2\le \e\Big\}\subset  \bigcup_{i, j} \hat D_m^{\e}(t_i,x_j)
\eeaa
Thus, by (\ref{Rm}),
\beaa
&&\hE\Big[\sup_{(t,x,\d, h)\in \hat D_N, |\d|+|h|^2\le \e}  |R_m(u,t,x,\d,h)|^p\Big] \le \sum_{i, j} \hE\Big[\sup_{(t,x,\d, h)\in \hat D_m^{\e}(t_i,x_j)}  |R_m(u,t,x,\d,h)|^p\Big]\\
&&\le C \sum_{i, j}  \e^{p(m+1)\over 2} \le C\e^{-1 - {d'(m+1)\over 2}} \e^{p(m+1)\over 2}  = C \e^{{(p-d')(m+1)\over 2} -1}.
\eeaa

\ms

(iv) Finally,  without loss of generality, we may assume ${p_* \over 1-\a} < p < p_0$. By (\ref{Rme}) we have
\beaa
&&\dbE\Big[\sup_{(t,x,\d,h) \in \hat D_N} \Big|{R_m(u,t,x,\d,h) \over (\d + |h|^2)^{m+\a\over 2}}\Big|^p\Big]\\
&=&\dbE\Big[\sup_{n\ge 0} \sup_{2^{-(n+1)}\le |\d|+|h|^2 \le 2^{-n}} \sup_{(t,x,\d,h) \in \hat D_N}\Big|{R_m(u,t,x,\d,h) \over (\d + |h|^2)^{m+\a\over 2}}\Big|^p\Big]\\
&\le& \sum_{n=0}^\infty \dbE\Big[\sup_{2^{-(n+1)}\le |\d|+|h|^2 \le 2^{-n}} \sup_{(t,x,\d,h) \in \hat D_N} \Big|{R_m(u,t,x,\d,h) \over (\d + |h|^2)^{m+\a\over 2}}\Big|^p\Big]\\
&\le& \sum_{n=0}^\infty2^{p(m+\a)(n+1)\over 2} \dbE\Big[\sup_{ |\d|+|h|^2 \le 2^{-n}} \sup_{(t,x,\d,h) \in \hat D_N} |R_m(u,t,x,\d,h)|^p \Big]\\
&\le& C \sum_{n=0}^\infty2^{p(m+\a)(n+1)\over 2}  2^{-n({(p-d')(m+1)\over 2} -1)} = C2^{p(m+\a)\over 2} \sum_{n=0}^\infty 2^{-{n\over 2}[p(1-\a)-p_*]}<\infty,
\eeaa
completing the proof.
\qed

\section{Extension to H\"{o}lder Continuous Case}
\setcounter{equation}{0}
\label{sect-Holder}

In this section we weaken the requirement $u\in \cH^{[m+2]}_{p_0}(\L)$ in Theorem \ref{thm-high} slightly, by replacing the highest order 
differentiability with a certain H\"{o}lder continuity. First recall the space $\cH^{[n]}_p(\L)$ defined in (\ref{cHn-space}). We shall now prove (\ref{inclusion}).
%and the fact 
% for any $\a \in (0,1)$, define 
% \bea
% \label{normKp2}
%\|u\|_{p,\a,T}&:=&  \|u\|_{p, T} + \hE\Big[\sup_{0\le t \le T, \d>0}{|u_{t+\d}  - u_t|^p \over \d^{p\a\over 2}}\Big]^{1/p},
% \eea 
%and,
% \bea
% \label{cHn2}
%&\cH^{[n]+\a}_p(\L) := \Big\{u\in \cH^{[n]}(\L): \|u\|_{n,p, \a,T} <\infty,  \forall T > 0\Big\}, ~~n\ge 0; \q \mbox{where}&\\
%&\dis\|u\|_{0,p,\a, T} := \|u\|_{p, \a, T},~ \|u\|_{1,p, \a, T} := \|u\|_{p,\a, T} +  \sum_{i=1}^{d} \|\pa_{\o^i}u\|_{p,\a, T}, &\nonumber\\
%&\dis\|u\|_{n,p, \a, T} := \|u\|_{p, \a, T} +\|\pa_t u\|_{n-2, p, \a, T} + \sum_{i=1}^{d} \|\pa_{\o^i}u\|_{n-1,p,\a, T}, \q n\ge 2.& \nonumber
%\eea
% It should be noted that $\cH^{[n+1]}_p(\L)$  is not a subspace of  $\cH^{[n]+\a}_p(\L)$. However, we have
\begin{lem}
\label{lem-inclusion}
Let $\a\in (0,1)$ and $p> {2\over 1-\a}$.  Then $\cH^{[n+2]}_p(\L) \subset \cH^{[n]+\a}_p(\L)$ for any $n\ge 0$.
\end{lem}

\no{\it Proof.} Without loss of generality, we shall only prove the case that $n=0$.  Let $u\in \cH^{[2]}_p(\L)$. First, for any $0\le t_0\le T$ and $0< \e < 1$,  applying functional It\^{o} formula (\ref{Ito}) and then Lemma \ref{lem-Iest},  we have
\beaa
\dbE\Big[ \sup_{0\le \d\le \e} \sup_{t_0\le t\le t_0+\e} |u_{t+\d} - u_t|^p\Big] \le C\e^{p\over 2}.
\eeaa
Then the lemma follows exactly the same arguments as in the proof of Theorem \ref{thm-high}.
\qed

In what follows we shall assume $u\in \cH^{[m+1]}_{p_0} (\L)\cap \cH^{[m]+\a_0}_{p_0}(\L)$ for some appropriate $m$, $p_0$ 
and $\a_0$. Note that in this case $R_m$ is still well-defined by (\ref{Taylorm}), however, we cannot use the representations (\ref{f-rep}) and (\ref{b-rep}) anymore because of their involvement of the $(m+2)$-th order derivatives.  

In order to estimate $R_m$ in this case, we first establish the following recursive representation. Recall that $\f_{s,t} := \f_t - \f_s$.

\begin{lem}
\label{lem-rep-Holder}  Let $u\in \cH^{[m+1]}_{2}(\L)$ for some $m\ge 1$. Then for any $\d>0$, it holds that:  
\bea
\label{f-rep1}
R_{m}(u, t, \d) &=& \sum_{|\th|\le m-2, \th_i\neq 0} \int_t^{t_\d} R_{m-2-|\th|} (\pa_t \cD^\th_\o u, t, s-t) \cI^\th_{s, t_\d} ds \\
&& + \sum_{|\th|=m-1, \th_i \neq 0} \int_t^{t_\d} \pa_t \cD^\th_\o u_s \cI^\th_{s, t_\d}  ds + \sum_{|\th|=m, \th_i \neq 0} \cI^\th_{t, t_\d}\big([\cD^\th_\o u]_{t,\cd}\big); \nonumber\\
\label{b-rep1}
R_{m}(u, t, -\d) &=& - \sum_{|\th|\le m} (-1)^{|\th|_0} R_{m-|\th|}(\cD^\th_\o u, t^-_\d, \d) \cI^{-\th}_{t^-_\d, t},\q t\ge \d.
\eea
\end{lem}

\no{\it Proof.}  (i)  We first prove (\ref{f-rep1}). We claim that, for $m\ge 2$,
\bea
\label{f-rep3}
R_{m}(u, t, \d) &=& \int_t^{t_\d} R_{m-2}(\pa_t u, t, s-t)ds + \sum_{i=1}^d\int_t^{t_\d} R_{m-1}(\pa_{\o^i} u, t, s-t) \circ dB^i_s.
\eea
Indeed, denote by $\tilde R_{m}(u, t, \d)$ the right side above, and notice the simple fact:
\bea
\label{fact6.1}
\sum_{|\th|\le m-2} \cD^\th_\o \pa_t u_t \int_t^{t_\d}  \cI^\th_{t,s}ds  +\sum_{|\th|\le m-1}  \cD^\th_\o\pa_{\o^i}u_t \sum_{i=1}^d \int_t^{t_\d} \cI^\th_{t,s} \circ dB^i_s=\sum_{1\le |\th|\le m} \cD^\th_\o u_t \cI^\th_{t,t_\d}. 
\eea
 Applying the functional It\^{o} formula (\ref{Ito2}) and by (\ref{Taylorm}), we have
\beaa
u_{t_\d} &=& u_t+\int_t^{t_\d} \pa_t u_s ds + \sum_{i=1}^d \int_t^{t_\d} \pa_{\o^i}u_s \circ dB^i_s\\
&=&u_t+ \int_t^{t_\d} \Big[\sum_{|\th|\le m-2} \cD^\th_\o \pa_t u_t \cI^\th_{t,s}  + R_{m-2}(\pa_t u, t, s-t)\Big] ds \\
&& + \sum_{i=1}^d \int_t^{t_\d} \Big[\sum_{|\th|\le m-1} \cD^\th_\o\pa_{\o^i}u_t \cI^\th_{t,s}  + R_{m-1}(\pa_{\o^i} u, t, s-t)\Big] \circ dB^i_s\\
&=&u_t+ \sum_{|\th|\le m-2} \cD^\th_\o \pa_t u_t \int_t^{t_\d}  \cI^\th_{t,s}ds  +\sum_{|\th|\le m-1}  \cD^\th_\o\pa_{\o^i}u_t \sum_{i=1}^d \int_t^{t_\d} \cI^\th_{t,s} \circ dB^i_s +  \tilde R_{m}(u, t, \d).
%\\&=&u_t +  \sum_{1\le |\th|\le m} \cD^\th_\o u_t \cI^\th_{t,t_\d}+  \tilde R_{m}(u, t, \d).
\eeaa
Then (\ref{f-rep3}) follows immediately from (\ref{fact6.1}) and (\ref{Taylorm}).

We now prove (\ref{f-rep1}). When $m=1$, the right side of (\ref{f-rep1}) becomes
\beaa
\int_t^{t_\d} \pa_t u_s ds + \sum_{i=1}^d \int_t^{t_\d} [\pa_{\o^i} u]_{t, s}\circ dB^i_s &=& \int_t^{t_\d} \pa_t u_s ds + \sum_{i=1}^d \int_t^{t_\d} \pa_{\o^i} u_s\circ dB^i_s - \sum_{i=1}^d \pa_{\o^i}u_t B^i_{t, t_\d}\\
&=& u_{t_\d} - u_t - \sum_{i=1}^d \pa_{\o^i}u_t B^i_{t, t_\d} = R_{1}(u, t, \d). 
\eeaa
For $m\ge 2$, applying (\ref{f-rep3}) repeatedly on the stochastic integral terms in (\ref{f-rep3}) we obtain
\bea
\label{f-rep4}
R_{m}(u, t, \d) &=&\sum_{|\th|\le m-2, \th_i\neq 0} \cI^\th_{t, t_\d}\Big( \int_t^\cd R_{m-2-\th}(\pa_t \cD^\th_\o u, t, s-t)ds\Big) \\
&&+ \sum_{|\th|=m-1, \th_i\neq 0} \cI^\th_{t, t_\d}\Big(R_1(\cD^\th_\o u, t, \cd-t)\Big).\nonumber
\eea
Applying stochastic Fubini theorem repeatedly, we have
\beaa
\cI^\th_{t, t_\d}\Big( \int_t^\cd R_{m-2-\th}(\pa_t \cD^\th_\o u, t, s-t)ds\Big) &=& \int_t^{t_\d} R_{m-2-\th}(\pa_t \cD^\th_\o u, t, s-t) \cI^\th_{s, t_\d}ds.
\eeaa
Moreover, note that $\cD^\th_\o u \in \cH^{[2]}_2(\L)$ for $|\th|=m-1$. Then
\beaa
&&R_1(\cD^\th_\o u, t, s-t) = [\cD^\th_\o u]_{t,s} - \sum_{j=1}^d \pa_{\o^j}\cD^\th_{\o} u_t B^j_{t, s}\\
&=& \int_t^s \pa_t\cD^\th_\o u_r dr + \sum_{j=1}^d \int_t^s \pa_{\o^j}\cD^\th_\o u_r \circ dB^j_r - \sum_{j=1}^d \pa_{\o^j}\cD^\th_{\o} u_t B^j_{t, s}\\
&=&\int_t^s \pa_t\cD^\th_\o u_r dr + \sum_{j=1}^d \int_t^s [\pa_{\o^j}\cD^\th_\o u]_{t,r} \circ dB^j_r.
\eeaa
Thus, (\ref{f-rep4}) leads to
\beaa
R_{m}(u, t, \d) &=&\sum_{|\th|\le m-2, \th_i\neq 0} \int_t^{t_\d} R_{m-2-\th}(\pa_t \cD^\th_\o u, t, s-t) \cI^\th_{s, t_\d}ds \\
&&+ \sum_{|\th|=m-1, \th_i\neq 0} \cI^\th_{t, t_\d}\Big(\int_t^\cd \pa_t\cD^\th_\o u_r dr + \sum_{j=1}^d \int_t^\cd [\pa_{\o^j}\cD^\th_\o u]_{t,r} \circ dB^j_r\Big),
\eeaa
which, together with  stochastic Fubini theorem again,  implies (\ref{f-rep1}) immediately.

\ms
(ii) We next prove (\ref{b-rep1}). By applying (\ref{Taylorm}) twice we have 
\beaa
&&R_{m}(u, t, -\d) = u_{t^-_\d} - \sum_{|\th|\le m} (-1)^{|\th|_0}\cD^\th_\o u_t \cI^{-\th}_{t^-_\d, t} \\
&=& u_{t^-_\d} - \sum_{|\th|\le m} (-1)^{|\th|_0} \Big[\sum_{|\tilde\th|\le m-|\th|} \cD^{\tilde\th}_\o \cD^\th_\o u_{t^-_\d} \cI^{\tilde\th}_{t^-_\d, t} + R_{m-|\th|}(\cD^\th_\o u, t^-_\d, \d)\Big] \cI^{-\th}_{t^-_\d, t} 
\eeaa
We now define
\bea
\label{Dm}
\D_m &:=& R_{m}(u, t, -\d)  + \sum_{|\th|\le m} (-1)^{|\th|_0} R_{m-|\th|}(\cD^\th_\o u, t^-_\d, \d) \cI^{-\th}_{t^-_\d, t} \\
&=&u_{t^-_\d} - \sum_{|\th|\le m} (-1)^{|\th|_0} \sum_{|\tilde\th|\le m-|\th|} \cD^{\tilde\th}_\o \cD^\th_\o u_{t^-_\d} \cI^{\tilde\th}_{t^-_\d, t}  \cI^{-\th}_{t^-_\d, t} .\nonumber
\eea
Denote $\hat \th := (\tilde\th, \th) = (\hat\th_1,\cds,\hat\th_n)$, then one can check that
\beaa
%\label{Dm2}
\D_m =  \sum_{1\le |\hat \th|\le m}  (-1)^{|\hat \th|_0+1}\cD^{\hat\th}_\o u_{t^-_\d}  a(\hat\th), &\mbox{where}& a(\hat\th) := \sum_{i=0}^n (-1)^i \cI^{(\hat\th_{1},\cds, \hat\th_i)}_{t^-_\d, t}  \cI^{(\hat\th_n,\cds, \hat\th_{i+1} )}_{t^-_\d, t}.
\eeaa
By setting $\f:=1$ in \reff{Rm-rep} we see that $a(\hat \th)=0$ for all  $1\le |\hat\th|\le m$.  Then $\D_m=0$ and we complete the proof.
\qed
\ms

With the above representations, we can now  extend Theorem \ref{thm-high}.
\begin{thm}
\label{thm-high-Holder} 
Assume that $u\in\cH^{[m+1]}_{p_0}(\L)\cap \cH^{[m]+\a_0}_{p_0}(\L) $ for some $m\ge 0$, $\a_0 \in (0, 1)$, and $p_0>  {2\over \a_0}$. The for any $p<p_0$ and $0< \a < \a_0-{2\over p_0}$,
it holds that, for any  $T>0$, 
\bea
\label{highRest-Holder}
\dbE\Big\{ \sup_{(t,\d) \in D^1_{[0,T]}}  \Big|{R_m(u, t; \d)\over |\d|^{m+\a\over 2}}\Big|^p\Big\} <\infty.
\eea
\end{thm}
\no{\it Proof.} Applying the representations (\ref{f-rep1}) and (\ref{b-rep1}), one  can easily prove by induction on $m$ that 
 \bea
  \label{temp-est-holder}
  \dbE\Big[\sup_{0<\d  \le \e}\sup_{t_0\le t\le t_0+\e} |R_m(u, t; \d )|^{p} + \sup_{0<\d  \le \e}\sup_{t_0\vee \d \le t\le t_0+\e} |R_m(u, t; -\d )|^{p}\Big]  \le C\e^{p(m+\a_0)\over 2},
  \eea
  extending (\ref{temp-est}). The result then follows from very similar arguments as those in \S\ref{sect-temp-proof}.  
  We leave the details to interested reader.
\qed

%Similarly,  for any $\a \in (0,1)$, define 
% \bea
% \label{normKp2-space}
%\|u\|_{p,\a,N}&:=&  \|u\|_{p, N} + \sup_{x\in K_N}\hE\Big[\sup_{0\le t < t'\le N}{|u(t,x) - u(t',x)|^p\over |t-t'|^{{p\a\over 2}} }\Big]^{1/p}\\
% &&+\sup_{0\le t\le N}\hE\Big[\sup_{x, x' \in K_N}{|u(t,x) - u(t,x')|^p\over|x-x'|^{p\a}} \Big]^{1/p},\nonumber
% \eea 
%and,
% \bea
% \label{cHn2-space}
%&\cH^{[n]+\a}_p (\hat\L):= \Big\{u\in \cH^{[n]}(\hat\L): \|u\|_{n,p, \a,N} <\infty,  \forall N > 0\Big\}, ~~n\ge 0; \q \mbox{where}&\\
%&\dis\|u\|_{0,p,\a, N} := \|u\|_{p, \a, N},~ \|u\|_{1,p, \a, N} := \|u\|_{p,\a, N} + \sum_{i=1}^{d'} \|\pa_{x_i} u\|_{p,\a, N} + \sum_{i=1}^{d} \|\pa_{\o^i}u\|_{p,\a, N}, &\nonumber\\
%&\dis\|u\|_{n,p, \a, N} := \|u\|_{p, \a, N} +\|\pa_t u\|_{n-2, p, \a, N} + \sum_{i=1}^{d'} \|\pa_{x_i} u\|_{n-1,p,\a, N} + \sum_{i=1}^{d} \|\pa_{\o^i}u\|_{n-1,p,\a, N}, & \nonumber
%\eea
%for $n\ge 2$ in the last line. 
% Applying Lemma \ref{lem-inclusion}, we have
%\bea
%\label{inclusion1}
%\mbox{$\cH^{[n+2]}_p(\hat\L) \subset \cH^{[n]+\a}_p(\hat\L)$ for any $n\ge 0$,  $\a\in (0,1)$ and $p> {2\over 1-\a}$.}
%\eea

Along the similar lines of arguments we can also weaken the assumptions of  Theorem \ref{thm-high-space}, the case for random field, to the H\"older conditions. Since the proof is a routine combination of the previous results, we omit it.
\begin{thm}
\label{thm-space-holder}
Assume that $u\in\cH^{[m+1]}_{p_0}(\L)\cap \cH^{[m]+\a_0}_{p_0}(\L) $ for some $m\ge 0$, $\a_0 \in (0, 1)$, and $p_0> p_*:= {1\over \a_0}[(m+\a_0)d' + 2]$. Then for any $p<p_0$ and $0< \a < \a_0[1-{p_*\over p_0}]$, (\ref{highRest-space}) holds  for any  $N>0$.
\end{thm}

\section{Application to (Forward) Stochastic PDEs}
\label{sect-SPDE}
\setcounter{equation}{0}

One of the main purposes of our study on the pathwise Taylor expansion is to lay the foundation for the notion of 
(pathwise) viscosity solution for stochastic PDEs and the associated forward path-dependent PDEs, which will be articulated
in our accompanying paper \cite{BMZ}. More precisely, we are particularly interested in 
the case when the random field $u$ is a (classical) solution (in standard sense) of the following SPDE:
\bea
\label{SPDE}
 d u(t,x)= f(t,x, u, \pa_x u, \pa^2_{x} u)dt+ g(t,x, u,\pa_x u)\circ dB_t,~t\ge 0,~\hP_0\mbox{-a.s.}
% u(t,x)=u_0(x)+\int_0^t f(s,x, u, \pa_x u, \pa^2_{x} u)ds+\int_0^t  g(s,x, u,\pa_x u)\circ dB_s,~t\ge 0,~\hP_0\mbox{-a.s.}
\eea
where $f(t,x, y, z, \g)$ and $g(t,x, y, z)$ are $\dbF$-progressively measurable and taking values in $\dbR$ and $\dbR^d$, respectively, with the variable $\o$ omitted as usual. 
Clearly, the SPDE (\ref{SPDE}) can be rewritten as   the following system of (forward) path dependent PDE (PPDE):
 \bea
 \label{PPDE}
 \pa_t u - f(t,x,u,\pa_x u, \pa_{xx}^2u) =0;\q \pa_\o u(t,x) = g(t,x,u,\pa_x u).
%\left\{\ba{lll}
 % \pa_t u - f(t,x,u,\pa_x u, \pa_{xx}^2u) =0;\\
 %  \pa_\o u(t,x) = g(t,x,u,\pa_x u);\\
 %   u(0,x) = u_0(x).
  %  \ea\right.
  \eea
  
%For simplicity, in this section we shall always assume
%\bea
%\label{cO}
%\cO := \dbR^{d'} &\mbox{and hence}& K_N := \{x\in \dbR^{d'} : |x|\le N\}.
%\eea
%Moreover, we define

%  \subsection{Pathwise Taylor expansion for SPDEs}

We will be particularly interested in the version of 
%TherIn order to establish our new notion of stochastic viscosity solution in \cite{BMZ}, we need to recast 
Theorem \ref{thm-high-space}, applyied to the solutions of (\ref{SPDE})  (or equivalently (\ref{PPDE})) in the case $m=2$.
% which will be the foundation of our notion of viscosity solution in
To this end, we first assume that $u$ is a  solution of SPDE (\ref{SPDE}) that is smooth enough in our sense. We shall also assume
that $g$ is sufficiently smooth. It is then clear that 
\bea
\label{pao}
\pa_\o u (t,x) = g(t,x, u(t,x), \pa_x u(t,x)).
\eea
Differentiating both sides above  in $x$ we get
(suppressing variables and noting that $g=(g_1,\cds, g_d)^T$): for $i=1,\cds, d'$, and $j=1,\cds, d$, 
\bea
\label{paxo}
\pa^2_{x_i\o^j} u = \pa_{x_i}\big(g_j(t,x, u , \pa_x u)\big)=\pa_{x_i} g_j+ \pa_y g_j \pa_{x_i} u + \sum_{k=1}^{d'}[\pa_{z_k} g_j] \pa^2_{x_ix_k} u,
\eea
or in matrix form:
\bea
\label{paxov}
\pa^2_{x\o} u &=& [\pa_{x_i\o^j}u]_{1\le i\le d', 1\le j\le d} \; =\;  \pa_{x} g + \pa_{x} u  [\pa_y g]^T+ \pa^2_{xx} u\; \pa_{z} g \in \dbR^{d'\times d}.
\eea
%where for $x\in\hR^{d'}$, $y\in \dbR^d$, $x \otimes y := xy^T= [x_iy_j]_{1\le i\le d', 1\le j\le d}$. 
Similarly, we can easily derive 
\bea
\label{paoo}
 \pa^2_{\o\o} u &=&  \pa_\o g + \pa_\o u [\pa_y g]^T+ [\pa^2_{x\o}u]^T[\pa_z g ]\\
 &=& \pa_\o g + g [\pa_y g]^T + \big[   \pa_{x} g + \pa_{x} u  [\pa_y g]^T+ \pa^2_{xx} u\; \pa_{z} g\big]^T[\pa_z g]. \nonumber
 \eea
In light of \reff{Taylorm-space} with $m=2$, we can now formally write down the pathwise Taylor expansion: 
 \bea
\label{Taylor2}
& &u(t+\d, x+h, \o)-u(t,x,\o) \nonumber\\
&=&\mbox{$\sum_{1\le|(\theta,\ell)|\le 2}\frac{1}{\ell !}(D_x^\ell D^\theta_\omega u)(t,x)h^\ell I^\theta_{t,t_\delta}+R(u,t,x,\delta,h)$}\nonumber\\
&=&f(t,x,u,\pa_x u, \pa^2_{xx} u)\d + \pa_x u \cd h +  g\cd \o^t_{t+\d}+ {1\over 2}\pa_{xx}^2 u : h h^T \nonumber\\
&& + \Big[\pa_{x} g + \pa_{x} u  [\pa_y g]^T+ \pa^2_{xx} u\; \pa_{z} g\Big]: h[B_{t,t+\d}]^T \\
&&+\Big[ \pa_\o g + g [\pa_y g]^T + \big[   \pa_{x} g + \pa_{x} u  [\pa_y g]^T+ \pa^2_{xx} u\; \pa_{z} g\big]^T[\pa_z g]\Big]: \underline B_{t, t+\d} + R(u, t,x, \d, h),\nonumber
\eea
for any $(t,x,\d,h) \in \hat D$, where the right hand side above is evaluated at $(t,x,\o)$, and $\underline B_{t, t+\d}$ is defined by (\ref{roughpath}). Applying Theorem \ref{thm-high-space} we then obtain the following result.
\begin{thm}
\label{thm-SPDE} 
Assume that SPDE (\ref{SPDE}) has a solution $u\in \cH^{[4]}_{p_0}(\hat\L)$ for some $p_0>  p_*:= 3 d'+2$, and that
the coefficient $g$ is regular enough so that all the derivatives in (\ref{paoo}) are well-defined. Let $R(u, t,x, \d, h)$ be determined (\ref{Taylor2}).
%by the following expansion:  
Then for any $p< p_0$ and $0<\a< 1-{p_*\over p_0}$,   the remainder $R$ satisfies  (\ref{highRest-space}) with $m=2$.
\end{thm}

  \begin{rem}
\label{rem-parabolicity}
{\rm The SPDE (\ref{SPDE}) can be written as the following It\^{o} form:
 \bea
\label{SPDE-Ito}
d u(t,x)= F(t,x, u, \pa_x u, \pa^2_{x} u)dt+ g(t,x, u,\pa_x u)\cd dB_t,~t\ge 0,~\hP_0\mbox{-a.s.}
\eea
where 
$$ F(t,x,y,z,\g) := f + {1\over 2} \tr\big(  \pa_\o g + g [\pa_y g]^T + \big[   \pa_{x} g + z  [\pa_y g]^T+ \g \; \pa_{z} g\big]^T[\pa_z g]\big).
$$
%\nonumber
%%u(t,x)=u_0(x)+\int_0^t F(s,x, u, \pa_x u, \pa^2_{x} u)ds+\int_0^t  g(s,x, u,\pa_x u)\cd dB_s,~~\dbP_0\mbox{-a.s.}
%\eea
It is thus natural to define the {\it  parabolicity} of the PPDE (\ref{PPDE}) as
\bea
\label{parabolic}
\pa_\g f = \pa_\g F- {1\over 2} [\pa_z g] [\pa_z g]^T \ge {\bf 0}.
\eea
Clearly, this is exactly the coercivity condition in the SPDE literature (see  e.g., Rozovskii \cite{Rozovski} and Ma-Yong \cite{MY1} in linear cases).
   \qed}
\end{rem}

%{\color{red}

We should note that the requirement  $u\in \cH^{[4]}_{p_0}(\hat\L)$ in Theorem \ref{thm-SPDE} is much stronger than $u$ 
being a classical solution, due to  the involvement of path derivatives. 
In the rest of this section we shall establish  the connection between the two concepts. To begin with let us recall 
some Sobolev spaces: for any $m\ge 0$ and $p\ge 1$,
\bea
\label{Hn}
W^m_p(\hat\L) := \{ u\in \dbL^0(\hat\L): \|u\|_{W^m_p,N} <\infty, \forall N>0\},\q W^m_*(\hat\L) := \cap_{p\ge 1} W^m_p(\hat\L),
\eea
where $\|u\|_{W^m_p,N}^p :=  \sup_{0\le t\le T} \sum_{|\ell|\le m}\hE\big[\int_{K_N}|\cD^\ell_x u(t,x)|^p dx\big]$. 

Next, we extend the  spatial derivatives slightly to those involving $(x,y,z, \g)$. To  this end, recall the multi-index set $\Th$  and the norm 
$|(\th,\ell)|$ on $\Th$ defined by (\ref{lnorm-space}). We define
%Since we are now considering spatial derivatives with respect to $(x, y,z, \g)$,  we extend the notation slightly:
%the following notation should be obvious:
\bea
\label{Dtildephi}
\cD^{(\th,\ell,\tilde\ell)}_{\o,(x,y,z)} \f :=\cD^\ell_x \cD^{\tilde\ell}_{(y,z)}\cD^\th_\o \f; ~\mbox{ and }~ 
\cD^{(\th,\ell,\hat\ell)}_{\o,(x,y,z, \g)} \f :=\cD^\ell_x \cD^{\hat\ell}_{(y,z, \g)}\cD^\th_\o \f,
\eea
%; ``length" of $\ell$ should be modified to that of  $(x,y,z,\g)$ accordingly. 
where, by a slight abuse of notation, $(\th,\ell,\tilde\ell)\in \Th\times \hN^{1+d'}$ and $(\th, \ell,\hat\ell)\in \Th\times \hN^{1+d'+d'\times d'}$, 
respectively; and $\f=\f(t,x,\o,y,z,\g)$ is any random field such that these derivatives exist.  Moreover, we define $|\tilde \ell|$ and $|\hat\ell|$ in an obvious way.
% Moreover, in order to  distinguish the $x$-derivatives and $(y,z)$- (and $(y,z,\g)$-)derivatives,  

We now state the main result of this section.
\begin{prop}
\label{prop-SPDE}
Let $m\ge 1$, $p> d'$, and denote $p_m := p(1+{m(m+1)\over 2})$. Assume that 

\ms

(i) for any $|\th|\le m-1$, $\cD^{(\th, \ell,\tilde \ell)}_{\o,(x,y,z)}g$ exists for all $|\ell|+|\tilde\ell| \le m-|\th|$, and is uniformly bounded when $|\tilde\ell|\ge 1$. Moreover,  $\cD^{\ell}_x \cD^\th_\o g$ is uniformly Lipschitz continuous in $(y,z)$  and $\cD^\ell_x\cD^\th_\o g(\cd, 0, {\bf 0}) \in W^{0}_{p_m}(\hat\L)$ for $|\ell|\le m-|\th|$.

\ms
(ii)  for any $|\th|\le m-2$, $\cD^{(\th,\ell,\hat \ell)}_{\o,(x,y,z,\g)}  f$ exists for all $|\ell|+|\hat\ell| \le m-1-|\th|$, and  is uniformly bounded when $|\hat\ell|\ge 1$. 
Moreover,  $\cD^\ell_x\cD^\th_\o f$ is uniformly Lipschitz continuous in $(y,z,\g)$ and $\cD^\ell_x\cD^\th_\o f(\cd, 0, {\bf 0}, {\bf 0}) \in W^0_{p_m}(\hat\L)$ for $|\ell|\le m-1-|\th|$. 

\ms
Let  $u \in W^{m+1}_{p_m}(\hat \L)$ be a classical solution (in the standard sense with differentiability in $x$ only) to SPDE (\ref{SPDE}),  %such that $\sum_{|\ell|\le 4} |\cD^{\ell}_x u| \in \dbL^{3p_0}(\hat\L)$, 
then  $u \in  \cH^{[m]}_{p}(\hat\L)$. 
%Moreover, for any $T>0$ we have
%\bea
%\label{uvanish}
%\lim_{N\to \infty} \sup_{0\le t\le T} \sup_{|x|\ge N} |u(t,x)|=0,\q \dbP_0\mbox{-a.s.}
%\eea
\qed
\end{prop}

To prove Proposition \ref{prop-SPDE} we need a technical lemma that would transform all  path derivatives to  the $x-$derivatives. 
Let us first introduce some
notations. For any random fields $\f=\f(t,x, \o, y,z,\g)$ and $u=u(t,x,\o)$, we define (suppressing variables)
\bea
\label{hatf}
\widehat \f (t,x,\o) := \f(t,x,\o, u, \pa_x u, \pa_{xx} u).
\eea
Next, for a given $(\bar \th, \bar\ell)\in \Th$, we define, with $(\th,\ell,\tilde\ell)\in \Th\times \hN^{1+d'}$, $(\th, \ell,\hat\ell)\in \Th\times \hN^{1+d'+d'\times d'}$, 
\bea
\label{SPDEA}
&& A_1:= A_1(\bar\th, \bar\ell):= \big\{ \ell: |\ell|\le |(\bar\th,\bar\ell)| \big\};\nonumber\\
&& A_2 :=A_2(\bar\th, \bar\ell):=\big\{(\th,\ell,\tilde\ell):  |\th|\le |\bar\th|-1,  |(\th,\ell,\tilde\ell)|\le |(\bar\th,\bar\ell)|-1\big\};\\
&&A_3 :=A_3(\bar\th, \bar\ell):=\big\{(\th,\ell,\hat\ell):  |\th|\le |\bar\th|-2,  |(\th,\ell,\hat\ell)|\le |(\bar\th,\bar\ell)|-2\big\}.\nonumber
\eea
%, in the above $\tilde \ell$ in $A_2$ (resp. $A_3$) denotes the order of $(y,z)$-(resp. $(y,z,\g)$-)derivative.

Let us now consider the random fields of the following form:
\bea
\label{term}
\psi:=\prod_{\ell\in A_1} (\cD^\ell_x u)^{a^1_\ell} \prod_{(\th,\ell,\tilde\ell)\in A_2} [\widehat{\cD^{(\th,\ell, \tilde\ell)}_{\o,(x,y,z)}g}]^{a^2_{(\th,\ell, \tilde\ell)}} \prod_{(\th,\ell,\hat\ell)\in A_3} [\widehat{\cD^{(\th,\ell,\hat\ell)}_{\o,(x,y,z,\g)}f}]^{a^3_{(\th,\ell, \hat\ell)}}
\eea
where,  $a^1_\ell, a^3_{(\th,\ell, \hat\ell)}\in \dbN$ and $a^2_{(\th, \ell, \tilde \ell)} \in \dbN^d$. We note that by definition the derivatives in (\ref{Dtildephi}) have the same dimension as the function $\f$. In particular, since 
$g$ is $\hR^d$-valued, the meaning of $[\h{\cD^{(\th,\ell,\tilde\ell)}g}]^a$, $a\in\hN^d$, in \reff{term} should be understood as 
that of $x^\ell$ defined in \reff{Dx}.

 Moreover, for each such $\psi$, we define its ``index", $\l(\psi)$,  by:
\bea
\label{SPDEl}
\l(\psi):= \sum_{\ell \in A_1} a^1_{\ell} + \sum_{(\th,\ell,\tilde\ell)\in A_2, \atop |\tilde\ell|=0} |a^2_{(\th,\ell, \tilde\ell)}| + 
\sum_{(\th,\ell,\hat\ell)\in A_3, \atop |\hat\ell|=0} a^3_{(\th,\ell, \hat\ell)}.
% \le 1+ |\bar\th|_0 |\bar\ell| + {|\bar\th|(|\bar\th|-1)\over 2}.
\eea
We remark that in the above we do not include the exponents when $|\tilde\ell|>0$ or $|\hat\ell|>0$, since the $(y,z,\g)$-derivatives are assumed to be bounded in Proposition \ref{prop-SPDE}. 

The following lemma will be crucial to the proof of Proposition \ref{prop-SPDE}. 
Since its proof is rather lengthy, we defer it to the end of the section in order not to disturb our discussion.
\begin{lem}
\label{lem-derivative}
Assume $f$ and $g$ are smooth enough with respect to all variables $(t,\o, x,y,z,\g)$, and $u$ is a classical solution (in standard sense) to SPDE \reff{SPDE} with sufficient regularity in $x$. Then, for any $(\bar\th, \bar\ell)\in \Th$,  $\cD^{\bar\ell}_x \cD^{\bar\th}_\o u$ is a linear combination of the terms in the form (\ref{term}). 

Moreover, for each term $\psi$, the following estimate holds for its index:
%we have the estimate for
\bea
\label{SPDElest}
\l(\psi)
%:= \sum_{\ell \in A_1} a^1_{\ell} + \sum_{(\th,\ell,\tilde\ell)\in A_2, \atop |\tilde\ell|=0} |a^2_{\th,\ell, \tilde\ell}| + 
%\sum_{(\th,\ell,\tilde\ell)\in A_3, \atop |\tilde\ell|=0} a^3_{\th,\ell, \tilde\ell} 
\le 1+ |\bar\th|_0 |\bar\ell| + {|\bar\th|(|\bar\th|-1)\over 2}.
\eea %$A_i = A_i(\bar\th, \bar\ell)$ are defined by
\end{lem}

\bs
Assuming this lemma we  now prove Proposition \ref{prop-SPDE}. 

\no[{\it Proof of Proposition \ref{prop-SPDE}.}] First, we recall so-called Morrey's inequality (cf. e.g., \cite{GT}) which states:
%. We state  it here for ready reference.
for any $\f: O\to \dbR$ that is in $W^1_p(O)$ (namely the generalized derivative $\pa_x \f$ is in $\dbL^p(O)$), where  $O\subset \dbR^{d'}$ is a 
bounded domain with $C^1$ boundary,  and any $p>d'$ and $0<\g < 1-{d'\over p}$, it holds that
\bea
\label{Morrey}
\sup_{x\in O} |\f(x)|^p + \sup_{x, x'\in O}\Big({ |\f(x)-\f(x')|\over |x-x'|^\g}\Big)^p \le C\int_O[|\f|^p+ |\pa_x \f|^p] dx.
\eea 
Now for $N>0$, recall the set $K_N$ defined by (\ref{KN}). Let $O$ be a domain with $C^1$ boundary such that $K_N \subset O \subset K_{N+1}$. 
From Morrey's inequality \reff{Morrey} we deduce that 
\beaa
\sup_{x\in K_N} |u(t,x)|^p\le C_N \int_{K_{N+1}} [|u(t,x)|^p+|\pa_x u(t,x)|^p] dx.
\eeaa
%where we used $K_{N+1}$ in the right side so that one 
Thus,  to prove the proposition it suffices to show that $\cD^\ell_x \cD^\th_\o u \in W^1_p(\hat\L)$ for all $|(\th, \ell)| \le m$, that is,
\bea
\label{SPDEclaim}
\cD^{\bar\ell}_x \cD^{\bar\th}_\o u \in W^0_p(\hat\L) &\mbox{for all}& |\bar\th| \le m,\q |(\bar\th, \bar\ell)|\le m+1.
\eea

To this end we fix $(\bar\th, \bar\ell)$ as in \reff{SPDEclaim}. If $m=1$,  then $\cD^{\bar\ell}_x \cD^{\bar\th}_\o u$ is either $\cD^{\bar\ell}_x u$ for $|\bar\ell| \le 2$, or $\cD^{\bar\ell}_x \pa_i u = \widehat{\cD^{\bar\ell}_x g^i}$ for some $i=1,\cds, d$ and $|\bar\ell| \le 1$, and one can check \reff{SPDEclaim} immediately. 

We thus assume $m\ge 2$. Denote, for $A_i := A_i(\bar\th, \bar\ell)$, $i=1,2,3$, we see that
\beaa
\xi :=\sum_{\ell \in A_1} |\cD^\ell_x u| +  \sum_{(\th,\ell,\tilde\ell)\in A_2, |\tilde\ell|=0} |\cD^\ell_x \cD^\th_\o g (\cd, 0, {\bf 0})| + \sum_{(\th,\ell,\hat\ell)\in A_3, |\hat\ell|=0} |\cD^\ell_x \cD^\th_\o f(\cd, 0, {\bf 0}, {\bf 0})|.
\eeaa
Note that,
\beaa
|\widehat{\cD^{\ell}_{x}\cD^{\th}_\o g}| &\le& | \cD^{\ell}_{x}\cD^{\th}_\o g(\cd, 0, 0)| + C[|u|+|\pa_x u|] \le C\xi,\q  (\th,\ell,\tilde\ell)\in A_2;\\
|\widehat{\cD^{\ell}_{x}\cD^{\th}_\o f}| &\le& |\cD^{\ell}_{x}\cD^{\th}_\o f (\cd, 0,0,0)| + C[|u|+|\pa_x u| + |\pa_{xx} u|]\le C\xi,\q  (\th,\ell,\hat\ell)\in A_3.
\eeaa
Note that $1+ |\bar\th|_0 |\bar\ell| + {|\bar\th|(|\bar\th|-1)\over 2} \le  1+{m(m+1)\over 2}$ for any $(\bar\th, \bar\ell)$ in \reff{SPDEclaim}. 
Applying Lemma \ref{lem-derivative}, one can then check that $|\cD^{\bar\ell}_x \cD^{\bar\th}_\o u| \le C\xi^{1+{m(m+1)\over 2}}$, which 
leads to \reff{SPDEclaim} immediately.
\qed

\ms
We now complete this section by proving  Lemma \ref{lem-derivative}.

\no[{\it Proof of Lemma \ref{lem-derivative}.}]  For simplicity, in this proof we assume $d=d'=1$. In particular, in this case 
$\bar \ell \in \dbN$ and thus $|\bar\ell| = \bar\ell$.  We first remark that if 
$\psi_1,\cds, \psi_n$ are the terms taking form of \reff{term}, then so is $\prod_{i=1}^n\psi_i$. Furthermore,  it holds that
\bea
\label{psiproperty}
  \l\Big(\prod_{i=1}^n\psi_i\Big) = \sum_{i=1}^n\l(\psi_i). 
\eea
We shall proceed in two steps.

{\it Step 1.} We first prove by induction on $|\bar\ell|$ that  $\cD^{\bar\ell}_x\big(\widehat f\big)$   is a linear combination of terms:
\bea
\label{fterm1}
\psi:=\prod_{|\ell|\le |\bar\ell|+2} (\cD^\ell_x u)^{a^1_\ell}  \prod_{|(\ell,\hat\ell)|\le |\bar\ell|} (\widehat{\cD^{(\ell,\hat\ell)}_{(x,y,z,\g)} f})^{a^3_{(\ell, \hat\ell)}};
\eea
and the index $\l(\psi)$ satisfies the estimate
\bea
\label{fterm2}
\l(\psi) := \sum_{|\ell|\le |\bar\ell|+2} a^1_\ell + \sum_{|\ell|\le |\bar\ell|, |\hat\ell|=0}  a^3_{(\ell, \hat\ell)} \le 1+|\bar\ell|.
\eea
Indeed,  when $|\bar\ell|=0$, we have $\psi := \widehat f$. Then all $a^1_\ell$'s and $a^3_{(\ell,\hat\ell)}$'s are equal to $0$ except $a^3_{(\ell,\hat \ell)} = 1$ for $|(\ell, \hat\ell)| = 0 $,  and thus $\l(\psi) = 1$.  

Assume the results hold true for $m$ and  $\bar\ell := \bar \ell' + 1$ with $|\bar\ell'|= m$.  Let $\psi'$ be a term in \reff{fterm1} corresponding to $\bar\ell'$. Then a typical term $\psi$ of $\cD^{\bar\ell} (\widehat f\big) = \pa_x \big[\cD^{\bar\ell'} (\widehat f\big)\big]$ should come from $\pa_x \psi'$. By \reff{psiproperty}, we now check the $x$-derivative of each factor of $\psi'$ and see its impact on  $\l$. 

First, for $|\ell|\le |\bar\ell'|+2$, we have $\pa_x  (\cD^\ell_x u) =  \cD^{\ell+1}_x u$, we see that $|\ell + 1|\le  |\bar\ell'|+2+2= |\bar\ell|+2$ and the corresponding $\l(\psi) = \l(\psi')$.  Next, for  $|(\ell',\hat\ell')|\le |\bar\ell'|$, 
\beaa
\pa_x \Big[\widehat{\cD^{(\ell',\hat\ell')}_{(x,y,z,\g)} f}\Big] &=& \widehat{[\pa_x\cD^{(\ell',\hat\ell')}_{(x,y,z,\g)} f]} + \widehat{[\pa_y\cD^{(\ell',\hat\ell')}_{(x,y,z,\g)}f]} \pa_x u\\
&& + \widehat{[\pa_z\cD^{(\ell',\hat\ell')}_{(x,y,z,\g)} f]} \pa_{xx} u + \widehat{[\pa_\g\cD^{(\ell',\hat\ell')}_{(x,y,z,\g)} f]} \pa_{xxx}u.
\eeaa
The derivatives of the $f$ terms are up to the order  $|(\ell', \hat\ell')|+1 \le  |\bar\ell'| +1 \le |\bar\ell|$, and those of  the $u$ terms are up to the order $3 \le |\bar\ell|+2$, so each term is still in the form of \reff{fterm1}. Moreover, the first three terms do not increase $\l$, while the last term increase $\l$ by $1$. Summarizing, we see that each term $\psi$ of $\pa_x \psi'$ is in the form of \reff{fterm1} and $\l(\psi) \le \l(\psi')+1$. Then we prove \reff{fterm2} for $|\bar\ell|$.

Similarly, we can prove that  $\cD^{\bar\ell}_x\big(\widehat g\big)$   is a linear combination of terms:
\bea
\label{gterm}
\psi:=\prod_{|\ell|\le |\bar\ell|+1} (\cD^\ell_x u)^{a^1_\ell}  \prod_{|(\ell,\tilde\ell)|\le |\bar\ell|} (\widehat{\cD^{(\ell,\tilde\ell)}_{(x,y,z)} g})^{a^2_{(\ell, \tilde\ell)}};
\eea
and the index satisfies the estimate: $ \l(\psi) := \sum_{|\ell|\le |\bar\ell|+1} a^1_\ell + \sum_{|\ell|\le |\bar\ell|, |\tilde\ell|=0}  a^2_{\ell, \tilde\ell} \le 1+|\bar\ell|$.

\ms

{\it Step 2.} We now prove the lemma  by induction on $|\bar \th|_0$. When $|\bar\th|_0=0$, the results are obvious. Assume the results hold true for $n$, and $\bar \th = (\th_1, \bar \th')$ with $|\bar\th'|_0=n$. Note that $\cD^{\bar\ell}_x \cD^{\bar\th}_\o u= \pa_{\th_1} \big(\cD^{\bar\ell}_x \cD^{\bar\th'}_\o u\big)$. Let $\psi'$ be a term in the form of \reff{term} corresponding to $(\bar\th', \bar\ell)$, then a typical term of $\cD^{\bar\ell}_x \cD^{\bar\th}_\o u$ should come from $\pa_{\th_1} \psi'$. We show that
\bea
\label{term-claim}
\mbox{$\psi$  is in the form of \reff{term}  corresponding to $(\bar\th, \bar\ell)$, and}~\l(\psi) \le \l(\psi') + |(\bar\th',\bar\ell)|.
\eea
This clearly implies  \reff{SPDElest} for $(\bar\th, \bar\ell)$.

We prove \reff{term-claim} in two cases. Denote $m:= |(\bar\th, \bar\ell)|$ and $m' := |(\bar\th', \bar\ell)|$. 

{\it Case 1.} $\th_1 = 0$. Then $|\bar\th|=|\bar\th'|+2$, $m=m'+2$ and  $\cD^{\bar\ell}_x \cD^{\bar\th}_\o u= \pa_t \big(\cD^{\bar\ell}_x \cD^{\bar\th'}_\o u\big)$.  By \reff{psiproperty}, we now check the $t$-derivative for each factor of $\psi'$ and see its impact on $\l$.

First, for $\ell \in A_1(\bar\th', \bar\ell)$, we have $\pa_t (\cD^{\ell}_x u) = \cD^{\ell}_x\Big( \widehat f\Big)$. Note that $|\ell|+2 \le m'+2 = m$ and $|\ell|\le m' =m-2$,  then \reff{fterm1}  implies that each term of $\cD^{\ell}_x\Big( \widehat f\Big)$ is in the form of \reff{term}.  Moreover, by  \reff{fterm2}, this differentiation increases the index $\l$ from $1$ up to $1+|\ell|\le 1+m'$. Then $\l(\psi) \le \l(\psi') + m'$.

Next, for $(\th,\ell,\tilde\ell)\in A_2(\bar\th', \bar\ell)$ with $m'\ge 1$ ($A_2(\bar\th', \bar\ell)$ is empty when $m'= 0$), we have 
\beaa
\pa_t \widehat{\big [\cD^{(\th,\ell,\tilde\ell)}_{\o,(x,y,z)} g\big]}
&=& \widehat{\big[\cD^{(\ell, \tilde\ell)}_{(x,y,z)}\pa_t\cD^{\th}_\o g\big]}  + \widehat{\big[\pa_y \cD^{(\ell, \tilde\ell)}_{(x,y,z)}\cD^{\th}_\o g\big]}\pa_t u + \widehat{\big[\pa_z \cD^{(\ell, \tilde\ell)}_{(x,y,z)}\cD^{\th}_\o g\big]}\pa_{xt} u \\
&=&\widehat{\big[\cD^{(\ell, \tilde\ell)}_{(x,y,z)}\pa_t\cD^{\th}_\o g\big]}  + \widehat{\big[\pa_y \cD^{(\ell, \tilde\ell)}_{(x,y,z)}\cD^{\th}_\o g\big]}\widehat f + \widehat{\big[\pa_z \cD^{(\ell, \tilde\ell)}_{(x,y,z)}\cD^{\th}_\o g\big]}\pa_{x} [\widehat f].
\eeaa
 The derivatives of $g$ are up to the order $|(\th,\ell,\tilde\ell)| + 2 \le m'-1+2 = m-1$, and its path derivatives are up to the order $|\th|+2 \le |\bar\th'|-1+2= |\bar\th|-1$,  then these terms  are in the form of \reff{term}. Moreover, since $m\ge 3$, by \reff{fterm1} one can easily see that all the terms of  $\widehat f$ and  $\pa_{x} [\widehat f]$ are in the form of \reff{term}. Furthermore, all the $g$-terms do not increase $\l$; the term $\widehat f$ increases $\l$ up to $1\le m'$. When $m'\ge 2$, the term $\pa_x[\widehat f]$ increase $\l$ up to $1+1\le m'$. When $m'=1$, we must have $|(\th,\ell,\tilde\ell)|=0$, then one can check straightforwardly that the $\l$ increases from $1$ to $2$, namely the increase is $1 =m'$. So in all the cases we have $\l(\psi) \le \l(\psi') + m'$. 

Finally,  for $(\th,\ell,\hat\ell)\in A_3(\bar\th', \bar\ell)$ with $m'\ge 2$ ($A_3(\bar\th', \bar\ell)$ is empty when $m'\le 1$), we have 
\beaa
\pa_t \widehat{\big [\cD^{(\th,\ell,\hat\ell)}_{\o,(x,y,z,\g)} f\big]}&=&\widehat{\big[\cD^{(\ell, \hat\ell)}_{(x,y,z,\g)}\pa_t\cD^{\th}_\o f\big]} + \widehat{\big[\pa_y\cD^{(\ell, \hat\ell)}_{(x,y,z,\g)}\cD^{\th}_\o f\big]} \widehat f\\
&& +\widehat{\big[\pa_z\cD^{(\ell, \hat\ell)}_{(x,y,z,\g)}\cD^{\th}_\o f\big]}\pa_{x} [\widehat f]+\widehat{\big[\pa_\g\cD^{(\ell, \hat\ell)}_{(x,y,z,\g)}\cD^{\th}_\o f\big]}\pa_{xx} [\widehat f].
\eeaa
The derivatives of $f$ are up to the order $|(\th,\ell)|  + 2 \le m'-2+2 = m-2$, and its path derivatives are up to the order $|\th|+2 \le |\bar\th'|-2+2= |\bar\th|-2$, then these terms  are in the form of \reff{term}.  Moreover, since $m\ge 4$, by \reff{fterm1} one can easily see that all the terms of  $\widehat f$, $\pa_{x} [\widehat f]$, and  $\pa_{xx} [\widehat f]$ are in the form of \reff{term}.  Furthermore, similarly to the $g$-case above, one can show that $\l(\psi) \le \l(\psi')+m'$.

{\it Case 2.} $\th_1 = 1$. Then $m=m'+1$ and  $\pa^{\bar\ell}_x \pa^{\bar\th}_\o u= \pa_\o \big(\pa^{\bar\ell}_x \pa^{\bar\th'}_\o u\big)$. By using \reff{gterm} and following similar arguments as in Case 1 we can easily prove the result.
 \qed

\section{Consistency with \cite{BBM}}
\label{sect-BBM}
\setcounter{equation}{0}

 In  this section we compare our stochastic Taylor expansions (\ref{Taylor2}) with those in our previous works 
 %It is worth noting that the   in 
 \cite{BM3, BBM} (in particular, the one in \cite{BBM}), and consequently unify them under the 
language of our path-derivatives. To be consistent with  \cite{BM3, BBM}, we assume in what follows that $d=1$, $\cO = \dbR^{d'}$, and that the coefficients $f$ and $g$ in (\ref{SPDE}) are {\it deterministic}. We should note that in this case we have
$A^t_{t+\d} =0$, and $\underline{B}_{t,t+\d}=\frac12 (\o^t_{t+\d})^2$.

% the respective
%coefficients in the corresponding Taylor expansions, in light of (\ref{paoo}), and  , assuming . 

We begin by recalling the definition of the  ``$n$-fold derivatives" introduced in \cite{BBM}.
\begin{defn}
\label{nfold} A random field $\zeta\in C^{0,n}(\hF^B,[0,T]\times
\dbR^{d'})$ is called ``$n$-fold" differentiable in the spatial variable
$x$ if there exist
% be a random field which is defined in an iterative way by a sequence of
$n$ random fields $\zeta_i\in C^{0,n}(\hF,[0,T]\times
\hR^{d'};\hR^{d_i})$, $2\leq i\leq n+1$, with $d_1=d'$ and $d_i\in\hN$,
$2\leq i\leq n+1$, and functions $F_i,G_i:[0,T]\times\hR^{d'}\times
\hR^{d_{i+1}}\rightarrow \hR^{d'}$, $i=1,\cds, n$,
% are assumed to satisfy the following regularity and growth conditions: \\
such that, denoting
%(for shortness of notation we put
$\zeta_1:=\zeta$, the following properties are satisfied:
\begin{itemize}
\item[(1)] $F_i,G_i\in C^\infty_{\ell,
p}$, $i=1,\cds,n$;

\item[(2)]  For $1\leq i\leq n$, it holds that
 \bea
 \label{zetai}
  \zeta_i(t,x)=\zeta_{i,0}(x)+\int_0^{t}F_{i}(s,x,\zeta_{i+1}(s,x))ds
 +\int_0^{t}G_{i}(s,x,\zeta_{i+1}(s,x))dB_s,
 \eea
for all $(t,x)\in [0,T]\times \hR^{d'}$, with $\zeta_{i,0}\in
C^2(\hR^{d'};\hR^{d_i})$, $1\leq i\leq n$.

\item[(3)] For any $1\leq i\leq n+1$, $N\geq 1$, and $|\ell|\le n$, 
%multi-index $j=(j_1,...,j_{d'})$ with $\vert j\vert=j_1+...+j_{d'}\leq n$,
 it holds that
$$\sup\{\cD^\ell_x\zeta_i(t,x)\vert, t\in
[0,T], \vert x\vert\leq N\}\in \cap_{p>1}L^p(\Omega, {\cF},\hP_0).
$$
\end{itemize}
We shall call $\z_i$, $i=2,\cds,n+1$ the ``generalized derivatives"
of $\z=\z_1$, with ``coefficients" $(F_i,G_i)$, $i=1,\cds,n$.
\end{defn}

The notion of the $n$-fold derivatives is particularly motivated by the structure of SPDE (\ref{SPDE}), or more precisely, (\ref{SPDE-Ito}). 
In fact, 
if we define $\zeta_1=u$, $\zeta_2=(u,\pa_x u,\pa_{xx}u)$, $F_1=F$, and $G_1=g$, then (\ref{zetai}) holds
for $i=1$. Moreover, if we assume that the coefficients $f$,  $g$ are sufficiently smooth so that the solution $u$ is 3-fold differentiable,
then all the coefficients of the 3-fold generalized derivatives can be determined by differentiating the equation (\ref{SPDE-Ito}) (in $x$)
repeatedly. 
Clearly, to  compare the Taylor expansion we need only compare the derivatives in Definition \ref{nfold} and the  path-derivatives defined in Definition \ref{Dw}.
For notational simplicity in what follows we shall assume  $d'=1$.

 To begin with, we note that  by Definition \ref{nfold} we have $\zeta_1=u$, $\zeta_2=(u,\pa_xu,\pa_{xx}u)$, $\zeta_3=(u,\pa_xu,\pa_{xx}u, \pa_{xxx}u, \pa_{xxxx}u)$. Thus
 we have $(F_1, G_1)=(F, g)$ and  the coefficients of the 2-fold derivative are $F_2=(F_2^{(1)},F_2^{(2)},F_2^{(3)})$ and $G_2=(G_2^{(1)},G_2^{(2)},G_2^{(3)})$, where $F_2^{(1)}=F$, $G_2^{(1)}=g$, and $(F^{(2)}_2, G^{(2)}_2)$ can be determined by
 differentiating the equation (\ref{SPDE-Ito}) with respect to $x$,  namely,
 \bea
\label{Dxu}
\pa_xu(t,x)=u'_0(x)+\int_0^t [F_2^{(2)}(t,x,\zeta_3(t,x))
%+F_yu_x+F_zu_{xx}+F_\g u_{xxx}
]ds+\int_0^t  G^{(2)}_2(t,x,\zeta_3(t,x))dB_s,
\eea
where, by direct calculation,  it is readily seen that 
\bea
\label{G22}
G^{(2)}_2(t,x,\zeta_3) =\pa_xg+\pa_yg\pa_xu+\pa_zg\pa^2_{xx}u.
\eea
Note that $g$ is deterministic, we have $\pa_\o g=0$ and, by (\ref{Dxu}), (\ref{G22}), and the definition of
path-derivative, 
%. So by Definition \ref{Dw}-(i) and (\ref{paoo})  we see that 
\bea
\label{Dxwg}
 \pa_{\o x}u=\pa_\o( \pa_xu)=G^{(2)}_2(t, x, \zeta_3)=\pa_xg+\pa_yg\pa_xu+\pa_zg\pa^2_{xx}u= \pa_{x}(\pa_\o u).
\eea
Now, applying functional  It\^o's formula (or ``chain rule") to $g(\cd, u,\pa_xu)$
%, and using 
%(\ref{SPDE})-(\ref{Dxwg}) 
we have
\bea
\label{gIto}
dg(t,x,u,\pa_xu)=F^g(t,x,\zeta_3)dt+[\pa_yg \pa_\o u+\pa_z g\pa_{\o x}u]dB_s,
\eea
where we simply denote the drift by $F^g$ as it is irrelevant to our argument. 
Note that, denoting $\bz=(y,z,\g)$, we can write 
$G_1(t,x,\bz)=G^{(1)}_2(t,x,\bz)=g(t,x,\bz)=g(t,x,y,z)$,
and $D_{\bz}g=(g_y,g_z,0)$.  Therefore
\bea
\label{DzG}
\pa_yg \pa_\o u +\pa_z g\pa_{\o x}u=\pa_y g g+\pa_zgG^{(2)}_2(t,x,\zeta_3)=\lan D_\bz g(t,x,\zeta_2), G_2(t,x,\zeta_3)\ran.
\eea
This, together with the fact $\pa_\o u=g$ and (\ref{gIto}),  shows that 
%$g$ is also an It\^o field, by Definition \ref{Dw} we have
\bea
\label{Dwwg}
\pa^2_{\o\o}u&=&\pa_\o g
= \lan D_\bz g, G_2(t,x,\zeta_3)\ran;\bs\\
\label{Dtu1}
\pa_tu&=&f(t,x,\zeta_2)-\frac12\pa^2_{\o\o}u=f(t,x,\zeta_2)-\frac12 \lan D_\bz g, G_2(t,x,\zeta_3)\ran.
\eea

We can now  recast the pathwise Taylor expansion of Buckdahn-Bulla-Ma \cite{BBM} in the new path-derivative language.
 Recall that $B^t_s(\o) :=\o_s - \o_{t}$, for $s\ge t$.
\begin{thm}
\label{StoTaylor}
Let $u$ be the classical solution to the SPDE (\ref{SPDE}) with deterministic coefficients
$f$ and $g$, and assume that it is 3-fold differentiable. 
Then, for every $\alpha \in
\big(\frac{1}{3},\frac{1}{2}\big)$ and  $m\in\hN$, there exist a
subset $\tilde{\O}_{\a,m}\subset\Omega$ with
$P\{\tilde\O_{\a,m}\}=1$, such that, for all $(t,y,\o)\in
[0,T]\times \overline{B}_m(0)\times\tilde\O_{\a,m}$,
the following Taylor expansion holds:
 \bea
 \label{taylorexp}
 u(t+h,y+k) -
 u(t,y)
 \neg\neg&\neg =\neg &\neg\neg ah + bB^t_{t+h} + \frac{c}{2} (B^t_{t+h})^2 + pk + \frac{1}{2} Xk^2 + qkB^t_{t+h}\nonumber\\
 && + (|h| + |k|^2)^{3\alpha}R_{\alpha,m}(t,t+h,y,y+k),
\eea
for all $(t+h, y+k)\in [0,T]\times \overline{B}_m(0)$. Here, the coefficients $(a, b, c, p, q, X)$ are given by
%with $(F_1,G_1)=(F,g)$, one has
\bea
\label{coeffi}
\left\{\ba{lll}
a=\pa_t u(t,x,\cd), \q b=\pa_\o u(t,x,\cd), \q c =\pa_{\o \o}u(t,x,\cd);\\
p=\pa_x u(t,x,\cd), \q q= \pa_{\o x}u(t,x,\cd), \q X= \pa_{xx}u(t,x,\cd).
\ea\right.
\eea
% such that $t+h\in [0,T]$ and all $k\in \dbR^d$ with $\vert k\vert\leq m$, where
Furthermore, the remainder of Taylor expansion $R_{\alpha,m}:
[0,T]^2\times (\hR^d)^2\times\Omega\mapsto \hR$ is a measurable
random field such that
 \begin{eqnarray}
 \label{Ram}
 \overline{R}_{\a,m}:=\sup_{t,s\in[0,T];~
 y,z\in\overline{B}_m(0)}|R_{\alpha,m}(t,s,y,z)|
%, \; 0 \leq t, s \leq T, \; |y|,|y+k|\leq m \}
 \in\cap_{p>1} L^{p}(\O, \cF, P).
 \end{eqnarray}
\end{thm}

{\it Proof.} Again, we assume $d'=1$. Since $u$ satisfies (\ref{SPDE}) and is 3-fold differentiable, by a direct application of Theorem 2.3 in \cite{BBM}, with $\zeta=\zeta_1=u$,
$\zeta_2=(u,u_x,u_{xx})$, and $F_1=F$, $G_1=g$, 
%f(t,x,u,u_x,u_{xx})=F_1(t,x,\zeta_2),\q
%g(t,x,u,u_x)=G_1(t,x,\zeta_2),
%\eea
we obtain a stochastic Taylor expansion (\ref{taylorexp}) with the following
coefficients
 \bea
 \label{taylors}
  a & = &f(t,x,\zeta_2)=F(t,x,\zeta_2) - \frac{1}{2} \lan D_\bz g(t,x,\zeta_2), G_2(t,x,\zeta_3)\ran,
  \nonumber\\
  b & = & g(t,x,\zeta_2), \q
  c  =  \lan D_\bz g(t,x,\zeta_2),G_2(t,x,\zeta_3)\ran, \nonumber\\
  p & = & \pa_x u(t,x), \q
  X  =  \pa_{xx}u(t,x), \\
  q & = & \pa_x g(t,x,\zeta_2) + \lan D_\bz g(t,x,\zeta_2)), D_x\zeta_2\ran. \nonumber
 \eea
Combining (\ref{pao}),  (\ref{Dxwg}), (\ref{Dwwg}) and (\ref{Dtu1}) we have
%recalling that $g$ does not depend on $\o$,
\beaa
a &=&f(t,x,\zeta_2)=F(t,x,\zeta_2) - {1\over 2}\lan D_\bz g(t,x,\zeta_2), G_2(t,x,\zeta_3)\ran= \pa_t u;\\
b &=& g (t,x,\zeta_2) = \pa_\o u;\\
c &=&\lan D_\bz g(t,x,\zeta_2), G_2(t,x,\zeta_3)\ran= \pa^2_{\o\o} u;\\
%p &=& \pa_x u(t,x);\q X= \pa^2_{xx} u(t.x);\\
q &=& \pa_xg(t,x,\zeta_2) +\lan D_\bz g(t,x,\zeta_2)), D_x\zeta_2\ran=\pa_x g + \pa_y g \pa_x u + \pa_z g \pa^2_{xx} u =   \pa^2_{x\o} u.
\eeaa
This proves  (\ref{coeffi}), whence the theorem.
\qed

  \end{document}